\documentclass[11pt,a4paper]{article}
\usepackage{amssymb}
\usepackage{rotating}
\usepackage{amsfonts}
\usepackage[colorlinks=true, urlcolor=blue, pdfborder={0 0 0}]{hyperref}
\usepackage{amsmath}
\usepackage{amsthm}
\usepackage{caption}
\usepackage{float}
\usepackage{graphicx}
\usepackage{wrapfig}
\usepackage{enumerate}
\usepackage{subcaption}
\usepackage{multicol}
\usepackage{srcltx}
\usepackage{array}
\newtheorem{theorem}{Theorem}[section]
\newtheorem{proposition}{Proposition}[section]
\newtheorem{lemma}{Lemma}[section]

\newtheorem{corollary}{Corollary}[section]

\newtheorem{remark}{Remark}
\newcommand{\wh}{\widehat}

\newcommand{\wt}{\widetilde}

\newcommand\cA{{\cal A}}

\newcommand\cG{{\cal G}}
\newcommand\cF{{\cal F}}
\newcommand\cL{{\cal L}}
\newcommand\cB{Q}
\newcommand\cP{{\cal P}}

\newcommand\cN{{\cal N}}

\newcommand\cR{{\cal R}}

\newcommand\cW{{\cal W}}

\newcommand\e{{\varepsilon}}

\def\bbr{{\mathbb R}}
\newcommand\Er{\mbox{Err}}
\newcommand{\aO}{\mbox{O}}
\newcommand{\ao}{\mbox{o}}

\def\text#1{\hbox{#1}}
\def\proof{{\noindent \bf Proof. }}

\def\endproof{\mbox{\ $\qed$}}

\def\E{{\bf E}}
\def\N{{\bf N}}

\def\P{{\bf P}}
\def\C{{\bf C}}

\def\p{{\bf p}}

\def\B{{\bf B}}

\def\U{{\bf U}}

\def\H{{\bf H}}
\def\R{{\bf R}}

\def\D{{\bf D}}

\def\r{{\bf r}}
\def\k{{\bf k}}

\def\u{{\bf u}}
\def\a{{\bf a}}
\def\s{{\bf s}}
\def\t{{\bf t}}
\def\b{{\bf b}}
\def\c{{\bf c}}
\def\l{{\bf l}}
\def\h{{\bf h}}
\newcommand\Trg{\mbox{Tr}}
\def\g{{\bf g}}
\def\v{{\bf v}}
\def\q{{\bf q}}

\def\m{{\bf m}}
\def\G{{\bf G}}
\def\Chi{{\bf 1}}
\def\d{\mathrm{d}}
\def\build #1_#2{\mathrel{\mathop{\kern 0pt #1}\limits_\zs{#2}}}
\newcommand{\zs}[1]{{\mathchoice{#1}{#1}{\lower.25ex\hbox{$\scriptstyle#1$}}
{\lower0.25ex\hbox{$\scriptscriptstyle#1$}}}}

\numberwithin{equation}{section}

\def\proof{{\noindent \bf Proof. }}
\def\endproof{\mbox{\ $\qed$}}

\usepackage[top=1.5in, bottom=1.5in, left=1.1in, right=1in]{geometry}

\begin{document}

\title{Adaptive
efficient robust sequential analysis for  autoregressive big data models
\thanks{
This research was supported by RSF, project no 20-61-47043
 (National Research Tomsk State University, Russia).}}

\author{Ouerdia Arkoun
\thanks{
Sup'Biotech, Laboratoire BIRL, 66 Rue Guy Moquet, 94800 Villejuif, France  and
Laboratoire de Math\'ematiques Raphael Salem, Normandie Universit\'e,
 UMR 6085 CNRS- Universit\'e de Rouen,  France,
  e-mail: ouerdia.arkoun@supbiotech.fr
}
,
Jean-Yves Brua
\thanks{
Laboratoire de Math\'ematiques Raphael Salem, Normandie Universit\'e,
 UMR 6085 CNRS- Universit\'e de Rouen,  France,
  e-mail :  jean-yves.brua@univ-rouen.fr
}
 and\\
 Serguei Pergamenshchikov\thanks{
 Laboratoire de Math\'ematiques Raphael Salem,
 UMR 6085 CNRS- Universit\'e de Rouen Normandie,  France
and
International Laboratory of Statistics of Stochastic Processes and
Quantitative Finance, National Research Tomsk State University,
 e-mail:
Serge.Pergamenshchikov@univ-rouen.fr } 
}


\maketitle

\begin{abstract}
In this paper we consider  high dimension models based on dependent observations defined through  autoregressive processes.  For  such models
we develop an adaptive efficient estimation method via
 the robust sequential model selection procedures.
 To this end, firstly,  using the Van Trees inequality, we obtain a sharp lower bound for robust risks in an explicit form given by the famous Pinsker constant (see in \cite{Pinsker1981, PchelintsevPergamenshchikovPovzun2021} for details). It should be noted, that for such models this constant is calculated for the first time.
Then, using the  weighted least square  method   and   sharp non asymptotic oracle inequalities from \cite{ArkounBruaPergamenchtchikov2019}
we provide the efficiency property in the minimax sense for the proposed estimation procedure, i.e.
we establish, that the upper bound  for its  risk 
  coincides with the obtained lower bound.  It should be emphasized that this  property  is obtained without using sparse conditions and in the adaptive setting when
 the  parameter dimension and model regularity are unknown. 
\end{abstract}



\newpage

\section{Introduction}

\subsection{Problem and motivations}

We study the observations model defined for $1\le j\le n$ through the following difference equation
\begin{equation}
\label{sec:In.1-hd-00}
y_\zs{j} = 
\left( 
\sum^{q}_\zs{i=1}\beta_\zs{i}\psi_\zs{i}(x_\zs{j})
\right)
y_\zs{j-1}+ \xi_\zs{j}
\,,
\quad
x_\zs{j}=a+\frac{(b-a)j}{n}\,,
 \end{equation}
where 
the initial value
  $y_\zs{0}$ is a non random known constant,
 $(\psi_\zs{i})_\zs{i\ge 1}$ are  known  linearly independent functions,
$a<b$ are fixed known constants
  and $(\xi_\zs{j})_\zs{j\ge 1}$ are i.i.d.
unobservable random variables
with an unknown density distribution  $\p$ from some functional class which will be specified later.

The problem
is to estimate the unknown parameters $(\beta_\zs{i})_\zs{1\le i\le q}$ in the high dimension setting, i.e. when the number of parameters
more than the number of observations, i.e. $q>n$.
Usually, in these cases  one uses one of two methods: Lasso algorithm 
or the Dantzig selector method (see, for example,  
\cite{HasttieFriedmanTibshirani2008} and the references therein).
It should be emphasized, that in the papers devoted to big data models the number of parameters $q$ must be known and,
moreover, usually it is assumed sparse conditions to provide  optimality properties.
Therefore, unfortunately, 
these methods can't be used  to estimate, for example, the 
number of parameters $q$. To overcome  these restrictions in this paper similar to the approach proposed in
 \cite{GaltchoukPergamenshchikov2021}
 we study this problem  in the nonparametric setting, i.e. we embed the observations 
 \eqref{sec:In.1-hd-00} in the general model defined as
\begin{equation}\label{sec:In.1}
y_\zs{j} = S(x_\zs{j}) y_\zs{j-1}+ \xi_\zs{j}\,,
\end{equation}
where $S(\cdot)\in\cL_\zs{2}[a,b]$ is unknown function. The nonparametric setting allows to consider
 the  models \eqref{sec:In.1-hd-00} with 
 unknown $q$ or even with $q=+\infty$. Note that the case when the number of parameters $q$ is unknown is one of 
  challenging problems
  in the signal and image processing theory 
(see, for example, \cite{BeltaiefChernoyarovPergamenchtchikov2020, BayisaZhouCronieYu2019}). 
So, now the problem is  to estimate the function $S(\cdot)$ 
on the basis of the observations \eqref{sec:In.1} under the 
condition that the noise  distribution $\p$ is unknown and belongs to some functional class 
$\cP$.
There is a number of papers which consider these  models (see, for example, \cite{Belitser2000, FanZhang2008,LuYangZhou2009}
and the references therein).
Firstly,  minimax estimation
problems
for the model \eqref{sec:In.1}
has been treated in 
\cite{ArkounPergamenshchikov2008, MoulinesPriouretRoueeff2005} in the nonadaptive case,
i.e. for the known regularity of the function
$S$.
Then, in  \cite{Arkoun2011}
and \cite{ArkounPergamenchtchikov2016}
 it is proposed to use the
sequential analysis 
for the adaptive  pointwise estimation  problem, i.e.
in the case when the  H\"older regularity is unknown.  Moreover, it turned out that only  the sequential methods  can provide the adaptive estimation for autoregressive models. 
That is why in this paper we use the adaptive sequential procedures
from \cite{ArkounBruaPergamenchtchikov2019}
 for the efficient estimation, which we study for the quadratic risks
 \begin{equation*}\label{sec:In.2-risk}
\cR_\zs{\p}(\wh{S}_\zs{n},S)\,=\,\E_\zs{\p,S}\|\wh{S}_\zs{n}-S\|^2\,,
\quad
\|S\|^2=\int^{b}_\zs{a}\,S^2(x)\d x\,,
\end{equation*}
where $\wh{S}_\zs{n}$ is an estimator of 
$S$ based on observations $(y_\zs{j})_\zs{1\le j\le n}$ and $\E_\zs{\p,S}$ is the expectation with respect to the distribution 
 law $\P_\zs{\p,S}$ of the process $(y_\zs{j})_\zs{1\le j\le n}$ given the distribution density $\p$ and the function $S$.
Moreover, taking into account that  $\p$ is unknown,  
 we use the robust risk defined  as
 \begin{equation}\label{sec:In.2-RUB-Risk}
\cR^{*}(\wh{S}_\zs{n},S)\,=\,
\sup_\zs{\p\in\cP}
\cR_\zs{\p}(\wh{S}_\zs{n},S)\,,
\end{equation}
where $\cP$ is a family of the distributions defined in Section \ref{sec:Main-Cnds}.

\subsection{Main tools}

To estimate the function $S$
in
model
\eqref{sec:In.1} we make use of the
model selection procedures proposed in 
\cite{ArkounBruaPergamenchtchikov2019}. These procedures are
based on the family of  the optimal poin-twise truncated sequential estimators 
from \cite{ArkounPergamenchtchikov2016} for which sharp oracle inequalities are obtained
 through the model selection method developed in \cite{GaltchoukPergamenshchikov2009a}.
In this paper, on the basis of these inequalities we show that the model selection procedures are efficient
in the adaptive setting for the robust quadratic risk \eqref{sec:In.2-RUB-Risk}. To this end, first of all,
we have to study the sharp lower bound for the these risks, i.e. we have to provide the best potential accuracy estimation
for the model \eqref{sec:In.1} which usually for quadratic risks is  given by the Pinsker constant 
(see, for example, in \cite{Pinsker1981, PchelintsevPergamenshchikovPovzun2021}). To do this we use the approach proposed in 
\cite{GaltchoukPergamenshchikov2009b,GaltchoukPergamenshchikov2011} based on the Van-Trees inequality. It turns out that for the model
\eqref{sec:In.1} the Pinsker constant
equals to the minimal quadratic risk value
for the filtration signal problem  studied in \cite{Pinsker1981}
multiplied by the coefficient obtained through the integration of the optimal pointwise estimation risk on the interval $[a,b]$.  This is a new result in the  nonparametric estimation theory
for the  statistical models with dependent observations.
Then, using the oracle inequality from \cite{ArkounBruaPergamenchtchikov2019}
and the weight least square estimation method we show that, for the 
model selection procedure  
the upper bound asymptotically coincides with the obtained Pinsker constant
without using the regularity properties of the unknown functions, i.e.
it is efficient in adaptive setting with respect to the robust risk
\eqref{sec:In.2-RUB-Risk}.

\subsection{Organization of the paper}

The paper is organized as follows. In  Section \ref{sec:Main-Cnds} we construct the sequential pointwise estimation procedures which allows us to pass from the autoregression model to the corresponding regression model, then 
 in Section \ref{sec:Ms} we develop the model selection method.
We announce the main results in Section \ref{sec:MainRes2}. In Section  \ref{sec:Siml}
we present Monte - Carlo results which numerically  illustrate the behavior of the proposed model selection procedures.  
In Section 
\ref{sec:VanTrees}
we show the Van-Trees inequality for the
model \eqref{sec:In.1}. We obtain the lower bound for the robust risk in Section \ref{sec:LowBnds-121} and
  in Section \ref{sec:MainResUp} we get the upper bound for the robust risk of the constructed sequential estimator. In Appendix  we give the all auxiliary  tools.

\section{Sequential point-wise estimation method}\label{sec:Main-Cnds}

To estimate the function $S$ in the model  \eqref{sec:In.1} on the interval $[a,b]$
we use the kernel sequentail estimators proposed in \cite{ArkounPergamenchtchikov2016}
 at the points $(z_\zs{l})_\zs{1\le l\le d}$
defined as
\begin{equation}\label{sec:Dt.00-00}
z_\zs{l}=a+\frac{l}{d}(b-a)\,,
\end{equation}
where $d$ is an integer valued function of $n$, i.e. $d=d_\zs{n}$, such that
$d_\zs{n}/\sqrt{n}\to 1$ as $n\to \infty$. Note that in this case the kernel estimator has the following form
\begin{equation*}
\label{sec:Kerne-Est-1}
\check{S}(z_\zs{l})=\frac{\sum^{n}_\zs{j=1}\,Q_\zs{l,j}\,y_\zs{j-1}\,y_\zs{j}}{\sum^{n}_\zs{j=1}\,Q_\zs{l,j}\,y^{2}_\zs{j-1}}
\quad\mbox{and}\quad
Q_\zs{l,j}=Q\left(\frac{x_\zs{j}-z_\zs{l}}{h}\right)
\,,
\end{equation*}
where $Q(\cdot)$ is  a kernel function and $h$ is a bandwidth. As is shown in \cite{ArkounPergamenshchikov2008}
to provide an efficient point wise estimation, the kernel function must be chosen as  the indicator  of the interval $]-1;1]$, i.e. 
$Q(u)=\Chi_\zs{]-1,1]}(u)$.  This means that we can rewrite the estimator \eqref{sec:Kerne-Est-1} as
\begin{equation}\label{sec:Sp.7}
\check{S}(z_\zs{l})=\frac{\sum^{k_\zs{2,l}}_\zs{j=k_\zs{1,l}}\,y_\zs{j-1}\,y_\zs{j}}{\sum^{k_\zs{2,l}}_\zs{j=k_\zs{1,l}}\,y^{2}_\zs{j-1}}
\,,
\end{equation}
where
$k_\zs{1,l}= [ n\wt{z}_\zs{l} - n\wt{h}]+1$ and
$k_\zs{2,l} =  [n\wt{z}_\zs{l} + n\wt{h}]\wedge n$, 
$[x]$ is the integer part of $x$,
$\wt{z}_\zs{l}=l/d$ and $\wt{h}=h/(b-a)$. 
To use the model selection method from \cite{ArkounBruaPergamenchtchikov2019}  
we need to obtain the uncorrelated stochastic terms 
in the kernel estimators for  the function $S$ at the points \eqref{sec:Dt.00-00}, i.e.
one needs to use the disjoint  observations  sets $(y_\zs{j})_\zs{k_\zs{1,l}\le j\le k_\zs{2,l}}$. To this end it suffices
to choose $h$ for which for all $2\le l\le d$ the bounds $k_\zs{2,l-1}< k_\zs{1,l}$, i.e. we set
\begin{equation}\label{sec:Sp.9}
h=\frac{b-a}{2 d}
\quad\mbox{and}\quad
\wt{h}=\frac{1}{2d}
\,.
\end{equation}

\noindent
Note that the main difficulty is that the kernel estimator is the non linear function of the observations due to the random  denominator.
To control this denominator we need to assume conditions providing the stability properties for the model \eqref{sec:In.1}.  To obtain the stable (uniformly with respect to the function $S$ ) model \eqref{sec:In.1}, we assume that for some fixed $0<\epsilon<1$ and $L>0$ the unknown function $S$ belongs to
the $\epsilon$-{\em stability set} introduced in \cite{ArkounPergamenchtchikov2016}
as
\begin{equation}\label{sec:Sp.1}
\Theta_\zs{\varepsilon,L} = \left\{S\in \C_\zs{1}([a,b],\bbr) : |S|_\zs{*} \le 1-\e
\quad\mbox{and}\quad
|\dot{S}|_\zs{*}\le L
 \right\}\,,
\end{equation}
where
$\C_\zs{1}[a,b]$ is the Banach space of continuously differentiable
$[a,b]\to\bbr$ functions and $|S|_\zs{*} = \sup_\zs{a\le x \leq b}|S(x)|$.
As is shown in 
\cite{ArkounBruaPergamenchtchikov2019}  $\forall S\in \Theta_\zs{\varepsilon,L}$  
$$
\sum^{k_\zs{2,l}}_\zs{j=k_\zs{1,l}}\,y^{2}_\zs{j-1}\approx  \frac{k_\zs{2,l}-k_\zs{1,l}}{1-S^{2}(z_\zs{l})}
\quad\mbox{as}\quad k_\zs{2,l}-k_\zs{1,l}\to\infty\,.
$$

\noindent
Therefore, to replace the denominator in \eqref{sec:Sp.7} with its limit we need firstly preliminary estimate the function $S(z_\zs{l})$. We estimate it as
\begin{equation}\label{sec:Sp.3}
\wh{S}_\zs{l}=\frac{\sum^{\iota_\zs{l}}_\zs{j=k_\zs{1,l}}\,y_\zs{j-1}\,y_\zs{j}}{\sum^{\iota_\zs{l}}_\zs{j=k_\zs{1,l}}\,y^{2}_\zs{j-1}}
\quad\mbox{and}\quad
\iota_\zs{l}=k_\zs{1,l}+\q\,,
\end{equation}
where $\q=\q_\zs{n}=[(n\wt{h})^{\mu_\zs{0}}]$
for some $0<\mu_\zs{0}<1$. Indeed, we can not use this estimator directly to replace the random denominator since in general it can be closed to one. By this reason we 
use its projection into the interval in  $]-1+\wt{\e},1-\wt{\e}[$, i.e. $\wt{S}_\zs{l} = \min (\max (\wh{S}_\zs{l},-1+\wt{\varepsilon}),1-\wt{\varepsilon})$ and
$\wt{\epsilon}=(2+\ln n)^{-1}$. 
Finally, omitting some technical details, we will replace the denominator \eqref{sec:Sp.7} with the threshold
  $H_\zs{l}$ defined  as

\begin{equation}\label{sec:Sp.8}
H_\zs{l}=\frac{1-\wt{\epsilon}}{1-\wt{S}^{2}_\zs{l}}
(k_\zs{2,l}-\iota_\zs{l})
\,.
\end{equation}
\noindent
It should be noted that $H_\zs{l}$ is a function the observations $y_\zs{k_\zs{1,l}},\ldots, y_\zs{\iota_\zs{l}}$. To replace the random denominator in \eqref{sec:Sp.7}
with the $H_\zs{l}$ we use the sequential estimation method through the following stopping time

\begin{equation}\label{sec:Sp.4}
 \tau_\zs{l} = \inf\{k>\iota_\zs{l}\,: \, \sum^{k}_\zs{j=\iota_\zs{l} +1} \u_\zs{j,l} \ge H_\zs{l} \}
 \,,
 \end{equation}
where $\u_\zs{j,l}=y^{2}_\zs{j-1}\Chi_\zs{\{\iota_\zs{l} +1\le j<k_\zs{2,l}\}}+H_\zs{l}\Chi_\zs{\{j=k_\zs{2,l}\}}$. It is clear that $\tau_\zs{l}\le k_\zs{2,l}$ a.s.
Now we define the sequential estimator as

\begin{equation}\label{sec:Sp.5}
S^{*}_\zs{l}=\frac{1}{H_\zs{l}}\,
\left(\sum^{\tau_\zs{l}-1}_\zs{j=\iota_\zs{l}+1}\,y_\zs{j-1}\,y_\zs{j}\,+\,\varkappa_\zs{l}\,y_\zs{\tau_\zs{l}-1}\,y_\zs{\tau_\zs{l}}\right)
\Chi_\zs{\Gamma_\zs{l}}\,,
\end{equation}
where
$\Gamma_\zs{l}=\{\tau_\zs{l} < k_\zs{2,l}\}$ and
 the correcting  coefficient $0<\varkappa_\zs{l}\le 1$ 
is defined as
\begin{equation}\label{sec:Sp.5-cr}
\sum^{\tau_\zs{l}-1}_\zs{j=\iota_\zs{l}+1}\u_\zs{j,l}
 + \varkappa^{2}_\zs{l}\u_\zs{\tau_\zs{l},l} = H_\zs{l}
\,.
\end{equation}

\noindent
To study robust properties of this sequential procedure  
similarly to
\cite{ArkounPergamenchtchikov2016}
we assume that  in the model \eqref{sec:In.1} the i.i.d. random
variables  $(\xi_\zs{j})_\zs{j\ge 1}$  have a density  $\p$
(with respect to the Lebesgue measure)  from the functional class $\cP$
defined as
\begin{align}\nonumber
\cP:=&
\left\{
\p\ge 0\,:\,\int^{+\infty}_\zs{-\infty}\,\p(x)\,\d x=1\,,\quad
\int^{+\infty}_\zs{-\infty}\,x\,\p(x)\,\d x= 0 \,,\right.
\\[2mm] \label{2.1}
&\quad\left.
\int^{+\infty}_\zs{-\infty}\,x^2\,\p(x)\,\d x= 1
\quad\mbox{and}\quad
\sup_\zs{l\ge 1}\,\frac{\int^{+\infty}_\zs{-\infty}\,|x|^{2l} \,\p(x)\,\d x}{l!\, \varsigma^{l}}
 \;\le 1
\,
\right\}\,,
\end{align}
where $\varsigma\ge 1$ is some fixed parameter, which may be a function
of the number observation $n$, i.e. $\varsigma=\varsigma(n)$, such that
 for any $\b>0$
\begin{equation}
\label{varsigma-cond-1-00}
\lim_\zs{n\to\infty}\,n^{-\b}\,\varsigma(n)
=0
\,.
\end{equation}

\noindent
  Note that the $(0,1)$  Gaussian density $\p_\zs{0}$
belongs to $\cP$.

\noindent
Now we can formulate
the following proposition from  \cite{ArkounBruaPergamenchtchikov2019} (Theorem 3.1).

\begin{proposition}\label{Pr.sec:Prs.stp.times.1}
For any $\b>0$ 
\begin{equation}\label{sec:Prs.stp.times.1}
\lim_\zs{n\to\infty}\,n^{\b}
\max_\zs{1\le l\le d}
\sup_\zs{S\in \Theta_\zs{\varepsilon,L}}\,\sup_\zs{\p\in \cP}
\P_\zs{\p,S}\left(\vert \wt{S}_\zs{l}-S(z_\zs{l})\vert>\epsilon_\zs{0}\right)
=0\,,
\end{equation}
where $\epsilon_\zs{0}=\epsilon_\zs{0}(n)\to 0$ as $n\to\infty$ such that 
$\lim_\zs{n\to\infty}n^{\gamma}\epsilon_\zs{0}=\infty$
for any $\gamma>0$.
\end{proposition}

\noindent 
Now we set
\begin{equation}\label{sec:In.df_est}
Y_\zs{l}=S_{\textcolor{red}{l}}^{*}\Chi_\zs{\Gamma}
\quad\mbox{and}\quad
\Gamma=
\cap^{d}_\zs{l=1}
\,\Gamma_\zs{l} 
\,.
\end{equation}

\noindent 
In Theorem 3.2 from  \cite{ArkounBruaPergamenchtchikov2019} 
it is shown that the probability of $\Gamma$ goes to zero uniformly more rapid than any power of the observations number  $n$, which is
formulated in the next proposition.

\begin{proposition}\label{Pr.sec:Prs.stp.times.22}
For any $\b>0$ the probability of the set $\Gamma$ satisfies the following
asymptotic equality
$$
\lim_\zs{n\to\infty}\,n^{\b}
\sup_\zs{S\in \Theta_\zs{\varepsilon,L}}
\P_\zs{\p,S}\left(\Gamma^{c}\right)
=0\,.
$$
\end{proposition}

\noindent
In view of this proposition  we can negligible the set $\Gamma^{c}$. So,
using the estimators \eqref{sec:In.df_est}
on the set $\Gamma$ we obtain the discrete time regression model

\begin{equation}\label{sec:In.1-11-1R}
Y_\zs{l}=S(z_\zs{l})+\zeta_\zs{l}
\quad\mbox{and}\quad
\zeta_\zs{l}=\eta_\zs{l}+\varpi_\zs{l}
\,,
\end{equation}
\noindent
in which
$$
\eta_\zs{l}
=\frac{
\sum^{\tau_\zs{l}-1}_\zs{j=\iota_\zs{l}+1}\u_\zs{j,l}\,\xi_\zs{j}+\,\varkappa_\zs{l}\u_\zs{\tau_\zs{l},l}\,\xi_\zs{\tau_\zs{l}}
}{H_\zs{l}}
\quad\mbox{and}\quad
\varpi_\zs{l}=\varpi_\zs{1,l}+\varpi_\zs{2,l}\,,
$$
where
$$
\varpi_\zs{1,l}=
\frac{
\sum^{\tau_\zs{l}-1}_\zs{j=\iota_\zs{l}+1}\u_\zs{j,l}\,
(S(x_\zs{j})-S(z_\zs{l}))
+
\varkappa^{2}_\zs{l}\u_\zs{\tau_\zs{l},l}\,
\,(S(x_\zs{\tau_\zs{l}})-S(z_\zs{\tau_\zs{l}}))
}{H_\zs{l}}
$$
and
$
\varpi_\zs{2,l}=
(\varkappa_\zs{l}-\varkappa^{2}_\zs{l})\u_\zs{\tau_\zs{l},l}\,S(x_\zs{\tau_\zs{l}})/H_\zs{l}$.
Note that 
the random variables 
$(\eta_\zs{j})_\zs{1\le j\le d}$
 (see Lemma A.2 in \cite{ArkounBruaPergamenchtchikov2019}), 
 for any $1\le l\le d$ and $\p\in\cP$ are such that
\begin{equation}\label{sec:def-eta-prs-25-03}
\E_\zs{\p,S}\left(\eta_\zs{l}\,\vert \cG_\zs{l}\right)=0
\,,
\quad
\E_\zs{\p,S}\left(\eta^{2}_\zs{l}\,\vert \cG_\zs{l}\right)=\sigma^{2}_\zs{l}
\quad\mbox{and}\quad
\E_\zs{\p,S}\left(\eta^{4}_\zs{l}\,\vert \cG_\zs{l}\right)\,\le \v^{*}\sigma^{4}_\zs{l}
\,,
\end{equation}
where  $\sigma_\zs{l}=H^{-1/2}_\zs{l}$, $\cG_\zs{l}=\sigma\{\eta_\zs{1},\ldots,\eta_\zs{l-1},\sigma_\zs{l}\}$ 
 and $\v^{*}$ is a fixed constant. 
 Note that
\begin{equation}
\label{sec:bound-sig-25}
\sigma_\zs{0,*}
\le
\min_\zs{1\le l\le d}\sigma^{2}_\zs{l}
\le 
\max_\zs{1\le l\le d}\sigma^{2}_\zs{l}
\le \sigma_\zs{1,*}
\,,
\end{equation}
where
$$
\sigma_\zs{0,*}=
\frac{1-\epsilon^{2}}{2(1-\wt{\epsilon})nh}
\quad\mbox{and}\quad
\sigma_\zs{1,*}=
\frac{1}{(1-\wt{\epsilon})(2nh-\q-3)}
\,.
$$


\begin{remark}
\label{Re.sec:seq-pr-**}
It should be summarized that we construct the sequential pointwise procedure 
\eqref{sec:Sp.4}  -- \eqref{sec:Sp.5} in two steps. First, we preliminary estimate the function 
$S(z_\zs{l})$ in \eqref{sec:Sp.3} on the observations $(y_\zs{j})_\zs{k_\zs{1,l}\le j\le \iota_\zs{l}}$ and through this estimator
we replace the random denominator in  \eqref{sec:Sp.5} with the threshold $H_\zs{l}$ in the second step when 
we construct
the estimation procedure on the basis of the observations $(y_\zs{j})_\zs{\iota_\zs{l}< j\le k_\zs{2,l}}$.  It should be noted also that 
 in the deviation \eqref{sec:In.1-11-1R} the main term $\eta_\zs{l}$ has a martingale form and the second one,
 as  it is shown in \cite{ArkounBruaPergamenchtchikov2019},
  is asymptotically small. It should be emphasized that namely these properties  allow us to develop effective estimation methods.   
\end{remark}

\section{Model selection}\label{sec:Ms}

Now to estimate the function $S$ we use the sequential model selection procedure from \cite{ArkounBruaPergamenchtchikov2019}
 for the regression
\eqref{sec:In.1-11-1R}. To this end, first  we choose the trigonometric basis  $(\phi_\zs{j})_\zs{j\ge\, 1}$ in $\cL_\zs{2}[a,b]$, i.e.
\begin{equation}\label{sec:In.5_trg}
\phi_\zs{1} = \frac{1}{\sqrt{b-a}}\,\,,\quad \phi_\zs{j}(x)= \sqrt{\frac{2}{b-a}}\, \Trg_\zs{j}\left(2\pi[j/2]\l_\zs{0}(x)\right)\,,\quad j\ge\,2\,,
\end{equation}
where the function $\Trg_\zs{j}(x)= \cos(x)$ for even $j$ and $\Trg_\zs{j}(x)= \sin(x)$ for odd $j$, 
and $\l_\zs{0}(x)=(x-a)/(b-a)$. Moreover, we choose
 the odd number $d$ of regression points \eqref{sec:Dt.00-00}, for example,
$d=2[\sqrt{n}/2]+1$. Then 
 the functions
$(\phi_\zs{j})_\zs{1\le j\le d}$ are orthonormal for the empirical inner product, i.e.
\begin{equation}\label{sec:Ms.4nn_1}
(\phi_i\,,\,\phi_\zs{j})_\zs{d}=
\frac{b-a}{d}\sum^{d}_\zs{l=1}\,\phi_i(z_\zs{l})\phi_\zs{j}(z_\zs{l})= \Chi_\zs{\{i=j\}}
\,.
\end{equation}

\noindent
It is clear that, the function $S$ can be represented as
\begin{equation}\label{sec:In.5_s--1}
S(z_\zs{l})=\sum^{d}_\zs{j=1}\,\theta_\zs{j,d}\,\phi_\zs{j}(z_\zs{l})
\quad\mbox{and}\quad
\theta_\zs{j,d}=
\left(S,\phi_\zs{j} \right)_\zs{d}
\,.
\end{equation}
We define the estimators for the coefficients $(\theta_\zs{j,d})_\zs{1\le j\le d}$ as
\begin{equation}\label{sec:In.5_est--1}
\wh{\theta}_\zs{j,d}=\frac{b-a}{d}
\sum^{d}_\zs{l=1}\,Y_\zs{l}\phi_\zs{j}(z_\zs{l})\,.
\end{equation}
From \eqref{sec:In.1-11-1R} we obtain immediately the following regression on the set $\Gamma$
\begin{equation}\label{sec:Ms.5}
\wh{\theta}_\zs{j,d}\,=\,\theta_\zs{j,d}\,+\,\zeta_\zs{j,d}
\quad\mbox{with}\quad
\zeta_\zs{j,d}=\sqrt{\frac{b-a}{d}}
\eta_\zs{j,d}+\varpi_\zs{j,d}\,,
\end{equation}
where
$$
\eta_\zs{j,d}=\sqrt{\frac{b-a}{d}}\sum^{d}_\zs{l=1}\,
\eta_\zs{l}\phi_\zs{j}(z_\zs{l})
\quad\mbox{and}\quad
\varpi_\zs{j,d}=\frac{b-a}{d}
\sum^{d}_\zs{l=1}\,\varpi_\zs{l}\,\phi_\zs{j}(z_\zs{l})\,.
$$
Through the Bounyakovskii-Cauchy-Schwarz we get that
\begin{equation}\label{sec:Ms.6}
|\varpi_\zs{j,d}|^{2}
\le\|\varpi\|^{2}_\zs{d}\,\|\phi_\zs{j}\|^{2}_\zs{d}=
\|\varpi\|^{2}_\zs{d}
\le (b-a)
 \frac{\varpi^{*}_\zs{n}}{n}
\,,
\end{equation}
where $\varpi^{*}_\zs{n}=\max_\zs{1\le l\le d}\varpi^{2}_\zs{l}$.
Note here, that as  it is shown in \cite{ArkounBruaPergamenchtchikov2019} (Theorem 3.3) for any $\b>0$ 
\begin{equation}\label{sec:Ms.6-lmt-01}
\lim_\zs{n\to\infty}\,
\frac{1}{n^{\b}}\,
\sup_\zs{\p\in\cP}\,
\sup_\zs{S\in\Theta_\zs{\epsilon,L}}
\E_\zs{\p,S}\varpi^{*}_\zs{n}\,\Chi_\zs{\Gamma}\,=\,0
\quad\mbox{for any}\quad
\b>0
\,.
\end{equation}

\noindent
To construct the model selection procedure we use 
 weighted least squares estimators 
 defined as
\begin{equation}\label{sec:Ms.8}
\wh{S}_\zs{\lambda}(t)=
\sum^{d}_\zs{l=1}\wh{S}_\zs{\lambda}(z_\zs{l})
\Chi_\zs{]z_\zs{l-1}, z_\zs{l}]}\,,\quad
\wh{S}_\zs{\lambda}(z_\zs{l})\,=\,\sum^{d}_\zs{j=1}\,\lambda(j)\,\wh{\theta}_\zs{j,d}\,\phi_\zs{j}(z_\zs{l})\,
\Chi_\zs{\Gamma}\,,
\end{equation}
where the weight vector $\lambda=(\lambda(1),\ldots,\lambda(d))'$ 
belongs to some finite set $\Lambda\subset [0,1]^{d}$, the prime denotes the transposition.
\noindent Denote by $\nu$  the cardinal number of the set $\Lambda$, for which we impose
the following condition.

\noindent $\H_\zs{1} :$
 {\em Assume that the number of the weight vectors $\nu$ 
as a function of $n$, i.e. $\nu=\nu_\zs{n}$, such that 
for any $\b>0$
the sequence $n^{-\b}\nu_\zs{n}\to 0$ as $n\to\infty$.
}


\noindent 
To choose a weight vector 
$\lambda\in\Lambda$ in \eqref{sec:Ms.8} we will use
 the following risk
\begin{equation}
\label{emperic_norm-110}
\Er_\zs{d}(\lambda)=\|\wh{S}_\zs{\lambda}-S\|^2_\zs{d}
=\frac{b-a}{d}\sum^{d}_\zs{l=1}(\wh{S}_\zs{\lambda}(z_\zs{l})-S(z_\zs{l}))^2\,.
\end{equation}

\noindent
Using \eqref{sec:In.5_s--1} and \eqref{sec:Ms.8} it can be represented as 
\begin{equation}\label{sec:Ms.10}
\Er_\zs{d}(\lambda)=
\sum^{d}_\zs{j=1}\,\lambda^2(j)\wh{\theta}^2_\zs{j,d}\,-
2\,\sum^{d}_\zs{j=1}\,\lambda(j)\wh{\theta}_\zs{j,d}\,\theta_\zs{j,d}\,+\,
\sum^{d}_\zs{j=1}\,\theta^2_\zs{j,d}\,.
\end{equation}
Since the coefficients $\theta_\zs{j,d}$ are unknown we can't minimize this risk directly
to obtain an optimal weight vector. To modify it  we set

\begin{equation}\label{sec:Ms.11}
\wt{\theta}_\zs{j,d}=
\wh{\theta}^2_\zs{j,d}-\frac{b-a}{d}\,\s_\zs{j,d}
\quad\mbox{with}\quad
\s_\zs{j,d}=\frac{b-a}{d}\,\sum^{d}_\zs{l=1}\,\sigma^2_\zs{l}\,\phi^2_\zs{j}(z_\zs{l})\,.
\end{equation}
Note here that in view of  
\eqref{sec:bound-sig-25}
-
\eqref{sec:Ms.4nn_1}
the last term can be estimated  as
\begin{equation}\label{sec:Ms.11++0}
\sigma_\zs{0,*}
\le
\s_\zs{j,d}\le \sigma_\zs{1,*}\,.
\end{equation}
Now, we modify the risk \eqref{sec:Ms.10} as
\begin{equation}\label{sec:Ms.12}
J_\zs{d}(\lambda)\,=\,\sum^{d}_\zs{j=1}\,\lambda^2(j)\wh{\theta}^2_\zs{j,d}\,-
2\,\sum^{d}_\zs{j=1}\,\lambda(j)\,\wt{\theta}_\zs{j,d}\,
+\,\delta P_\zs{d}(\lambda)\,,
\end{equation}
where the coefficient $0< \delta< 1$  will be chosen later
and
 the penalty term is 
\begin{equation}\label{sec:Pen-term}
P_\zs{d}(\lambda)=
\frac{b-a}{d}
\sum^{d}_\zs{j=1}\lambda^2(j)\s_\zs{j,d}\,.
\end{equation}
\noindent
Now using  \eqref{sec:Ms.12}
we define the sequential model selection procedure as
\begin{equation}\label{sec:Ms.13}
\wh{\lambda}=\mbox{argmin}_\zs{\lambda\in\Lambda}\,J_\zs{d}(\lambda)
\quad\mbox{and}\quad  
\wh{S}_\zs{*}=\wh{S}_\zs{\wh{\lambda}}\,.
\end{equation}

\noindent
To study the efficiency property
we specify the weight coefficients
$(\lambda(j))_\zs{1\le j\le n}$ as it is proposed, for example, in \cite{GaltchoukPergamenshchikov2009b}.
First, for some $0<\varepsilon<1$
introduce the two dimensional grid
to adapt to  the unknown parameters (regularity and  size)  of the Sobolev ball, i.e. we set
\begin{equation}\label{sec:Ga.0}
\cA=\{1,\ldots,k^*\}\times\{\varepsilon,\ldots,\m\varepsilon\}\,,
\end{equation}
where $\m=[1/\varepsilon^2]$. We assume that both
parameters $k^*\ge 1$ and $\varepsilon$ are functions of $n$, i.e.
$k^*=k^*(n)$ and $\varepsilon=\varepsilon(n)$, such that
\begin{equation}\label{sec:Ga.1}
\left\{
\begin{array}{ll}
&\lim_\zs{n\to\infty}\,k^*(n)=+\infty\,,
\quad\lim_\zs{n\to\infty}\,\dfrac{k^*(n)}{\ln n}=0\,,\\[6mm]
&
\lim_\zs{n\to\infty}\,\varepsilon(n)=0
\quad\mbox{and}\quad
\lim_\zs{n\to\infty}\,n^{\b}\varepsilon(n)\,=+\infty
\end{array}
\right.
\end{equation}
for any $\b>0$. One can take, for example, for $n\ge 2$
\begin{equation}\label{sec:Ga.1-00}
\varepsilon(n)=\frac{1}{ \ln n }
\quad\mbox{and}\quad
k^*(n)=k^{*}_\zs{0}+ [\sqrt{\ln n}]\,,
\end{equation}
where $k^{*}_\zs{0}\ge 0$ is some fixed integer number.
 For each $\alpha=(k, \t)\in\cA$, we introduce the weight
sequence
$
\lambda_\zs{\alpha}=(\lambda_\zs{\alpha}(j))_\zs{1\le j\le d}$
with the elements
\begin{equation}\label{sec:Ga.3}
\lambda_\zs{\alpha}(j)=\Chi_\zs{\{1\le j<j_\zs{*}\}}+
\left(1-(j/\omega_\alpha)^{k}\right)\,
\Chi_\zs{\{ j_\zs{*} \le j\le \omega_\zs{\alpha}\}}\,,
\end{equation}
 and
  $$
 \omega_\zs{\alpha}=
 \omega_\zs{*}
 +
 (b-a)^{2k/(2k+1)}
 \left(
 \frac{(k+1)(2k+1)}{\pi^{2k}k}
 \t\,n
 \right)^{1/(2k+1)}
 \,.
 $$
 \noindent
 Here, $j_\zs{*}$ and $\omega_\zs{*}$ are such that $j_\zs{*}\to\infty$,
  $j_\zs{*}=\ao\left((n/\varepsilon)^{1/(2k+1)}\right)$ and $\omega_\zs{*}=\aO(j_\zs{*})$ as $n\to\infty$.
In this case we set
$\Lambda\,=\,\{\lambda_\zs{\alpha}\,,\,\alpha\in\cA\}$\,.
Note, that these weight coefficients are
  used in \cite{KonevPergamenshchikov2012, KonevPergamenshchikov2015}
 for continuous time regression
models  to show the asymptotic efficiency.
It will be noted that in this case the cardinal of the set $\Lambda$ is  
$\nu=k^{*} m$. It is clear that the properties \eqref{sec:Ga.1}
imply the condition $\H_\zs{1}$.
In \cite{ArkounBruaPergamenchtchikov2019} we showed the following result.

\begin{theorem}\label{Th.sec:OrIn.1}
Assume that the conditions \eqref{varsigma-cond-1-00}
and
$\H_\zs{1}$
hold.
 Then for any
 $n\ge\, 3$, any
$S\in\Theta_\zs{\varepsilon,L}$ and any $ 0 <\delta \leq 1/12$, 
 the procedure \eqref{sec:Ms.13} with the coefficients  \eqref{sec:Ga.3}
satisfies the following oracle inequality
\begin{equation}
\label{sec:Mrs.1--00-12}
\cR^{*}(\wh{S}_*,S)\leq\frac{(1+4\delta)(1+\delta)^2}{1-6\delta} \min_\zs{\lambda\in\Lambda} 
\cR^{*}(\wh{S}_\lambda,S)
+
\,
\frac{\B^{*}_\zs{n}}{\delta n}
\,,
\end{equation}
where the term $\B^{*}_\zs{n}$ is such that 
$\lim_\zs{n\to\infty}\, n^{-\b}\B^{*}_\zs{n}=0$ for any $\b>0$.
\end{theorem}

\begin{remark}
\label{Re.sec:ModSelec-00+}
In this paper we will use the inequality   \eqref{sec:Mrs.1--00-12} to study  efficiency properties for the model selection procedure 
\eqref{sec:Ms.13}  with the weight coefficients  
\eqref{sec:Ga.3} in adaptive setting, i.e. in the case when the regularity of the function $S$ \eqref{sec:In.1} is unknown.
\end{remark}

\section{Main results}\label{sec:MainRes2}
First, to study the minimax properties for the estimation problem for the model \eqref{sec:In.1}
 we need to introduce some functional class. To this end
 for any fixed
$\r>0$ and $\k\ge 2$
we set
\begin{equation}\label{2.5}
\cW_\zs{\k,\r}=\left\{f\in\Theta_\zs{\varepsilon,L}:\,
\sum_\zs{j=1}^{+\infty}\,\a_\zs{j}\,\theta^{2}_\zs{j}\le \r\right\}\,,
 \end{equation}
where
$\a_\zs{j}=\sum^{\k}_\zs{l=0}\,(2\pi[j/2]/(b-a))^{2l}$,
  $(\theta_\zs{j})_\zs{j\ge 1}$ are the trigonometric Fourier coefficients in $\cL_\zs{2}[a,b]$, i.e.
$
\theta_\zs{j}=(f,\phi_\zs{j})=\int^{b}_\zs{a} f(x)\phi_\zs{j}(x)\d x$
and $(\phi_\zs{j})_\zs{j\ge 1}$ is the trigonometric basis defined in \eqref{sec:In.5_trg}.
  It is clear that we can represent this functional class as the Sobolev ball
$$
\cW_\zs{\k,\r}=\left\{f\in\Theta_\zs{\varepsilon,L}:\,
\sum_\zs{j=0}^{\k}\,\|f^{(j)}\|^2\le \r
\right\}\,.
$$

\noindent
Now, for this set we define the normalizing coefficients
\begin{equation*}\label{2.7-0}
\l_\zs{\k}(\r)=
\l_\zs{\k}(\r)\,=\,((1+2\k)r)^{1/(2\k+1)}\,
\left(\frac{\k}{\pi (\k+1)}\right)^{2\k/(2\k+1)}
\end{equation*}
and
\begin{equation}\label{2.7}
\varsigma_\zs{*}=
\varsigma_\zs{*}(S)=\int_\zs{a}^{b} (1-S^{2}(u)) \d u
\,.
\end{equation}
It is well known that
in regression models with the functions
  $S\in \cW_\zs{\k,\r}$ 
 the minimax  convergence rate  is 
$n^{-2\k/(2\k+1)}$ (see, for example,  \cite{GaltchoukPergamenshchikov2009b, KonevPergamenshchikov2009}
and the references therein). Our goal in this paper is to show the same property for the non parametric auto-regressive models
\eqref{sec:In.1}. First we have to obtain a lower bound for the risk \eqref{sec:In.2-RUB-Risk}   over all possible
estimators $\Xi_\zs{n}$, i.e. any  measurable function of the observations  $(y_\zs{1},\ldots,y_\zs{n})$.  


\begin{theorem}\label{Th.2.1-1}
For the model \eqref{sec:In.1} 
the robust risk  \eqref{sec:In.2-RUB-Risk} normalized by the coefficient $\upsilon(S)=((b-a)\varsigma)^{-2\k/(2\k+1)}_\zs{*}$
can be estimated from below as
 \begin{equation}\label{2.9}
\liminf_\zs{n\to\infty}\inf_{\wh{S}_\zs{n}\in\Xi_\zs{n}}\,n^{2\k/(2\k+1)}\sup_\zs{S\in \cW_\zs{\k,\r}}
\,\upsilon(S)\cR^{*}(\wh{S}_\zs{n},S)\,
\ge \l_\zs{\k}(\r)\,.
 \end{equation}
\end{theorem}


\noindent
Now to study the procedure \eqref{sec:Ms.13} we have to add some condition on the penalty coefficient $\delta$ which provides sufficiently small
 penalty term in \eqref{sec:Ms.12}.

\noindent
$\H_\zs{2}$ : {\em Assume that the parameter $\delta$ is a function of $n$, i.e. $\delta=\delta_\zs{n}$ such that 
$\lim_\zs{n\to\infty}\delta_\zs{n}=0$ and
$\lim_\zs{n\to\infty}\,n^{-\b}\delta_\zs{n}
=0$
for any $\b>0$.
}

\begin{theorem}\label{Th.2.2-2} 
Assume that the conditions $\H_\zs{1}$ -- $\H_\zs{2}$ hold. 
Then the model selection procedure
  $\wh{S}_\zs{*}$ defined in 
 \eqref{sec:Ms.13} with the weight vectors \eqref{sec:Ga.3}
 admits the following asymptotic upper bound
 \begin{equation*}\label{2.8}
\limsup_\zs{n\to\infty}\,n^{2\k/(2\k+1)}
\sup_\zs{S\in \cW_\zs{\k,\r}}\,\upsilon(S)\,
\cR(\wh{S}_\zs{*},S)\,
\le \l_\zs{\k}(\r) \,.
 \end{equation*}
\end{theorem}

\noindent
Now
Theorems 
\ref{Th.2.1-1} - \ref{Th.2.2-2} imply the following efficiency property.

\begin{corollary}\label{Co.2.2-2} 
Assume that the conditions $\H_\zs{1}$ -- $\H_\zs{2}$ hold. 
The model selection procedure
  $\wh{S}_\zs{*}$ defined in 
 \eqref{sec:Ms.13} and \eqref{sec:Ga.3}
is efficient, i.e.
 \begin{equation}\label{Co.2.8}
 \lim_\zs{n\to\infty}
 \frac{
 \inf_{\wh{S}_\zs{n}\in\Xi_\zs{n}}\,\sup_\zs{S\in \cW_\zs{\k,\r}}\,\upsilon(S)
\,\cR^{*}(\wh{S}_\zs{n},S)}{ 
\sup_\zs{S\in \cW_\zs{\k,\r}}\,\upsilon(S)\cR(\wh{S}_\zs{*},S)}
 =1
 \,.
 \end{equation}
\noindent
 Moreover,
 \begin{equation}\label{Co.2.8--1}
 \lim_\zs{n\to\infty}\,n^{2\k/(2\k+1)}
\sup_\zs{S\in \cW_\zs{\k,\r}}\,\upsilon(S)\,
\cR^{*}(\wh{S}_\zs{*},S)\,
= \l_\zs{\k}(\r) \,.
\end{equation}
\end{corollary}

\begin{remark}\label{Re.2.2}
Note that the limit equalities \eqref{Co.2.8} and 
\eqref{Co.2.8--1} imply that the function 
$\l_\zs{\k}(\r)/\upsilon(S)$
is the
minimal value of the normalized asymptotic quadratic robust risk, i.e.
the Pinsker constant 
in this case. We remind that the coefficient $\l_\zs{\k}(\r)$
is the well known Pinsker constant for the "signal+standard white noise" model obtained in \cite{Pinsker1981}. 
Therefore,
the Pinsker constant for the model \eqref{sec:In.1}
is represented by the Pinsker constant for the 
"signal+white noise" model  in which the noise intensity is given by
the function 
\eqref{2.7}.
\end{remark}

Now we assume that in the model  \eqref{sec:In.1} the functions $(\psi_\zs{i})_\zs{i\ge 1}$ are
orthonormal in $\cL_\zs{2}[a,b]$, i.e. 
$
(\psi_\zs{i},\psi_\zs{j})
=\Chi_\zs{\{i=j\}}$.
We use the estimators \eqref{sec:Ms.8}
to estimate the parameters $\beta=(\beta_\zs{i})_\zs{i\ge 1}$ as
$
\wh{\beta}_\zs{\lambda}=(\wh{\beta}_\zs{\lambda,i})_\zs{i\ge 1}$
and
$\wh{\beta}_\zs{\lambda,i}=(\psi_\zs{i},\wh{S}_\zs{\lambda})$.
Then, similarly we use the selection model procedure \eqref{sec:Ms.13}
as
\begin{equation}
\label{selmd-1}
\wh{\beta}_\zs{*}=(\wh{\beta}_\zs{*,i})_\zs{i\ge 1}
\quad\mbox{and}\quad
\wh{\beta}_\zs{*,i}=(\psi_\zs{i},\wh{S}_\zs{*})
\,.
\end{equation}

\noindent
It is clear, that in this case
$
\vert \wh{\beta}_\zs{\lambda}-\beta \vert^{2}=
\sum^{\infty}_\zs{i=1}(\wh{\beta}_\zs{\lambda,i}-\beta_\zs{i})^2
=\Vert \wh{S}_\zs{\lambda}-S\Vert^{2}$
and
$\vert \wh{\beta}_\zs{*}-\beta \vert^{2}
=\Vert \wh{S}_\zs{*}-S\Vert^{2}$.

\par
Note, that Theorem \ref{Th.sec:OrIn.1} implies that the estimator \eqref{selmd-1}
is optimal in the sharp oracle inequality sense which is established in the following theorem.
\begin{theorem}\label{Th.sec:bgd-1.1}
For any 
$S\in\Theta_\zs{\varepsilon,L}$,
 $n\ge\, 3$ and $ 0 <\delta \leq 1/12$, 
$$
\cR^{*}(\wh{\beta}_\zs{*},\beta)
\le\,
\frac{(1+4\delta)(1+\delta)^2}{1-6\delta}\,
 \min_\zs{\lambda\in\Lambda}
\cR^{*}(\wh{\beta}_\zs{\lambda},\beta)
 +\frac{\B^{*}_\zs{n}}{\delta n}\,,
$$
where 
$\cR^{*}(\wh{\beta},\beta)=\sup_\zs{\p\in\cP}\E_\zs{\p,S}\vert \wh{\beta}-\beta \vert^{2}$ and
  $\B^{*}_\zs{n}$ satisfies the limit property mentioned in Theorem \ref{Th.sec:OrIn.1}.
\end{theorem}

\par
Note now, that
Theorems
\ref{Th.2.1-1} and  \ref{Th.2.2-2}
imply the efficiency property for the estimate \eqref{selmd-1} based on
the model selection procedure
 \eqref{sec:Ms.13}
constructed  with 
the penalty threshold $\delta$ satisfying the
 condition $\H_\zs{2}$.
\begin{theorem}\label{Th.sec: bgd121.1-EFF}
Then the estimate \eqref{selmd-1}
is asymptotically efficient, i.e.
$$
\lim_\zs{n\to\infty}\,n^{2\k/(2\k+1)}\,
\sup_\zs{S\in \cW_\zs{\k,\r}}\,
\upsilon(S)
\cR^{*}(\wh{\beta}_\zs{*},\beta)
= l_\zs{k}(\r)
$$
and
\begin{equation}\label{sec:Ae.5--10}
\lim_\zs{n\to\infty}\frac{\inf_\zs{\wh{\beta}_\zs{n}\in\Xi_\zs{n}}
\,
\sup_\zs{S\in \cW_\zs{\k,\r}}\,\upsilon(S)
\cR^{*}(\wh{\beta}_\zs{n},\beta)
}
{\sup_\zs{S\in \cW_\zs{\k,\r}}\,
\upsilon(S)
\cR^{*}(\wh{\beta}_\zs{*},\beta)
}
=1\,,
\end{equation}
where  $\Xi_\zs{n}$ is the set  estimators for  $\beta$ based on the observations $(y_\zs{j})_\zs{1\le j\le n}$.
\end{theorem}



\begin{remark}\label{Re.2.3}
It should be noted that we obtained the efficiency property \eqref{sec:Ae.5--10}  for the big data autoregressive model  \eqref{sec:In.1}
without using the parameter dimension $q$ or sparse conditions usually used for such models (see, for example, in 
\cite{HasttieFriedmanTibshirani2008}). 
\end{remark}

\section{Monte - Carlo simulations}\label{sec:Siml}

In this section we present the numeric  results obtained through the Python soft  
for the model \eqref{sec:In.1} in which $(\xi_\zs{j})_\zs{1\le j\le n}$ are i.i.d. $\cN(0,1)$ random variables and  $0\le x\le 1$, i.e. $a=0$ and $b=1$.
In this case we simulate the  
 model selection procedure \eqref{sec:Ms.13} with the weights \eqref{sec:Ga.3} in which 
$k^*= 150+\sqrt{\ln n}$, $\m=[\ln^{2} n]$,  $\varepsilon=1/\ln n$. Moreover, the parameters $j_\zs{*}$ and
 $\omega_\zs{*}$ are chosen as
 $$
 j_\zs{*}=\frac{\underline{\omega}}{200+\ln \underline{\omega}}\,,
 \qquad
 \underline{\omega}
 =
 \ln n
 +
 \left(
 \frac{(k+1)(2k+1)}{\pi^{2k}k}
 \t\,n
 \right)^{1/(2k+1)}
$$
and $\omega_\zs{*}=j_\zs{*}+\ln n$.
 First we study the model \eqref{sec:In.1} with
$S_\zs{1}(x)=0,5 \cos(2\pi x)$ and then for the function
$$
S_\zs{2}(x)=0,1+\sum^{q}_\zs{j=1}\,\frac{\cos(2\pi j x)}{(j+3)^{2}}
\quad\mbox{and}\quad
q=100000
\,.
$$
In the model selection procedures we use
$d=2[\sqrt{n}/2]+1$ points in \eqref{sec:Dt.00-00}.

\newpage

%
%

\begin{figure}[h!]\label{fig.n-200}
\begin{subfigure}[b]{0.5\textwidth}

  \includegraphics[width=\columnwidth,height=6cm]{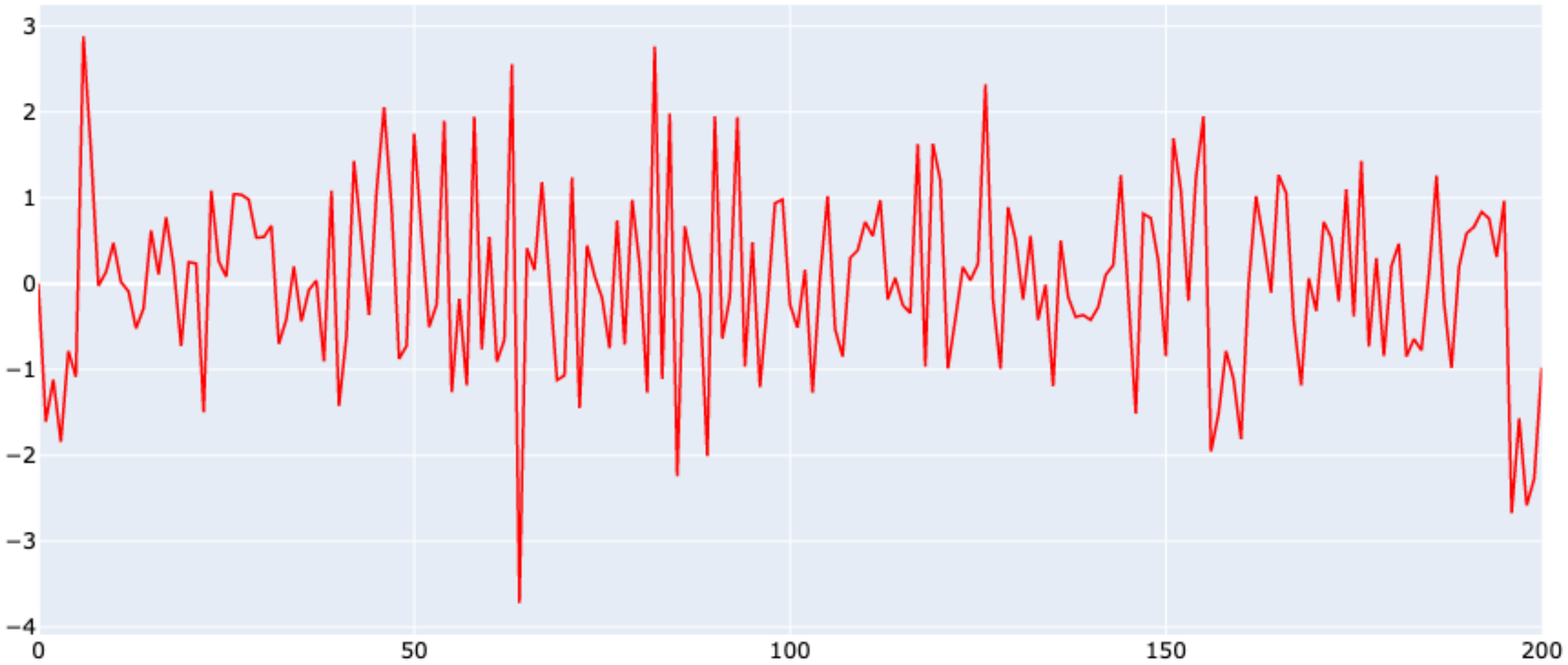}
\caption{\small Observations for $n=200$}%
\label{fig.1-n-200}
\end{subfigure}%
\begin{subfigure}[b]{0.5\textwidth}
\includegraphics[width=\columnwidth,height=6cm]{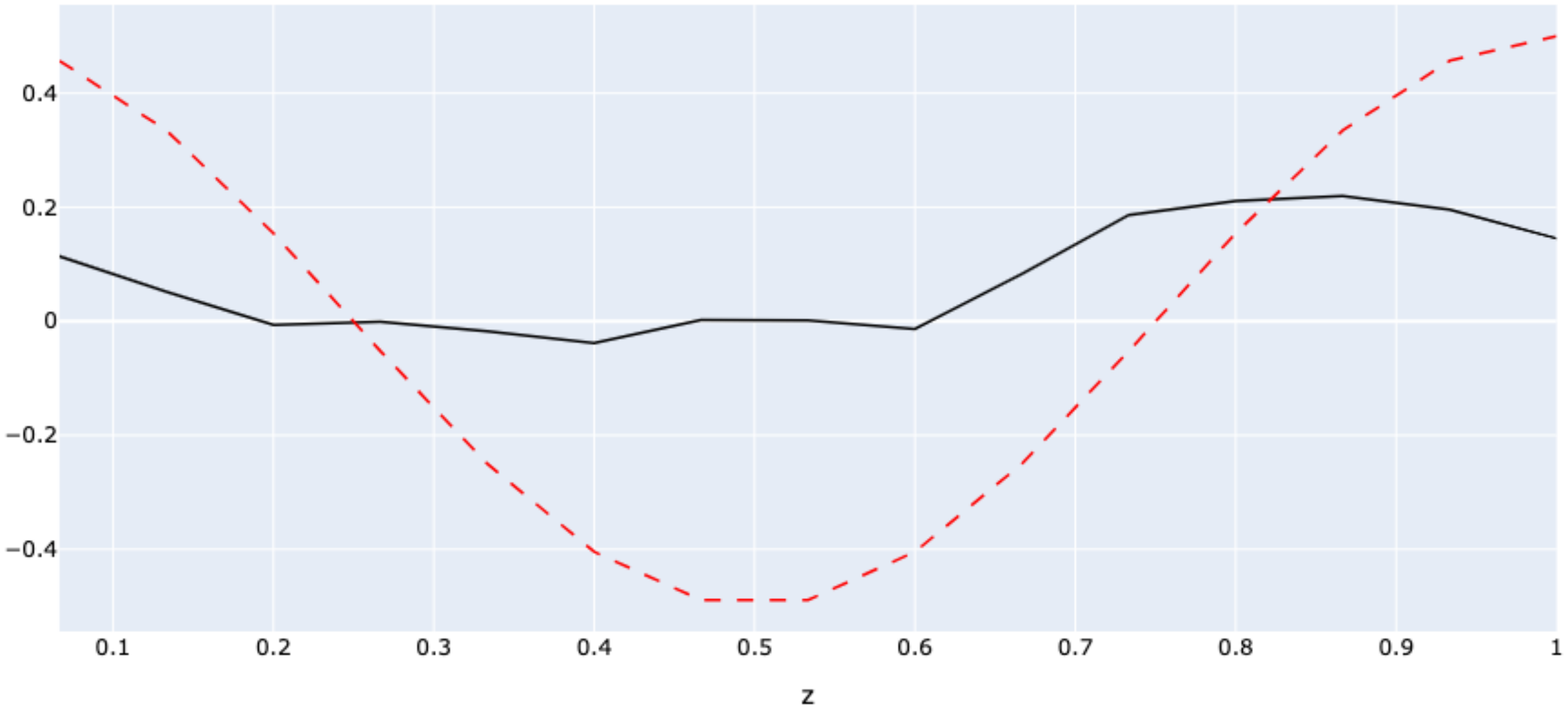}
\caption{\small Estimator of $S_\zs{1}$ for $n=200$}%
\label{fig.2-n-200}
\end{subfigure}
\caption{\small Model selection for $n=200$}\label{fig.n-200} 
\end{figure}

\begin{figure}[h]\label{fig.n-500}
\begin{subfigure}[b]{0.5\textwidth}

  \includegraphics[width=\columnwidth,height=6cm]{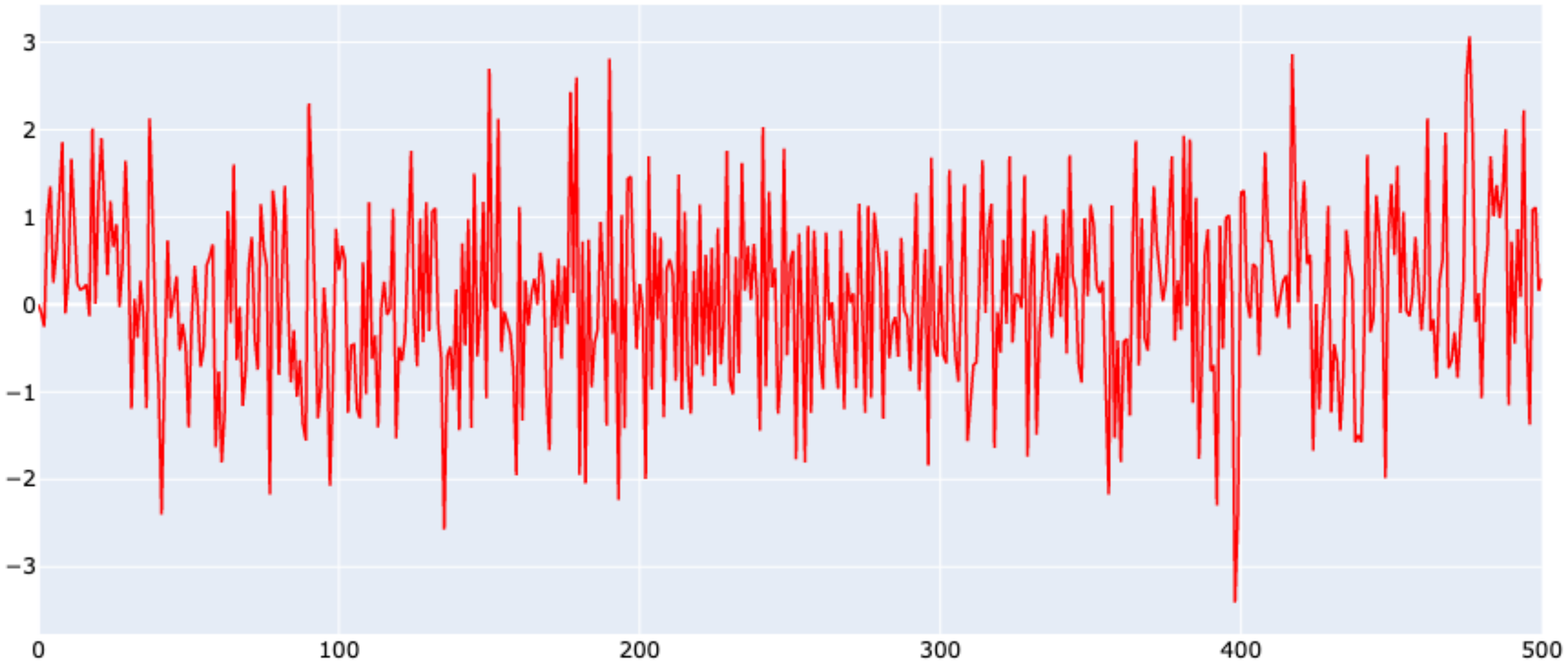}
\caption{\small Observations for $n=500$}%
\label{fig.1-n-500}
\end{subfigure}%
\begin{subfigure}[b]{0.5\textwidth}
\includegraphics[width=\columnwidth,height=6cm]{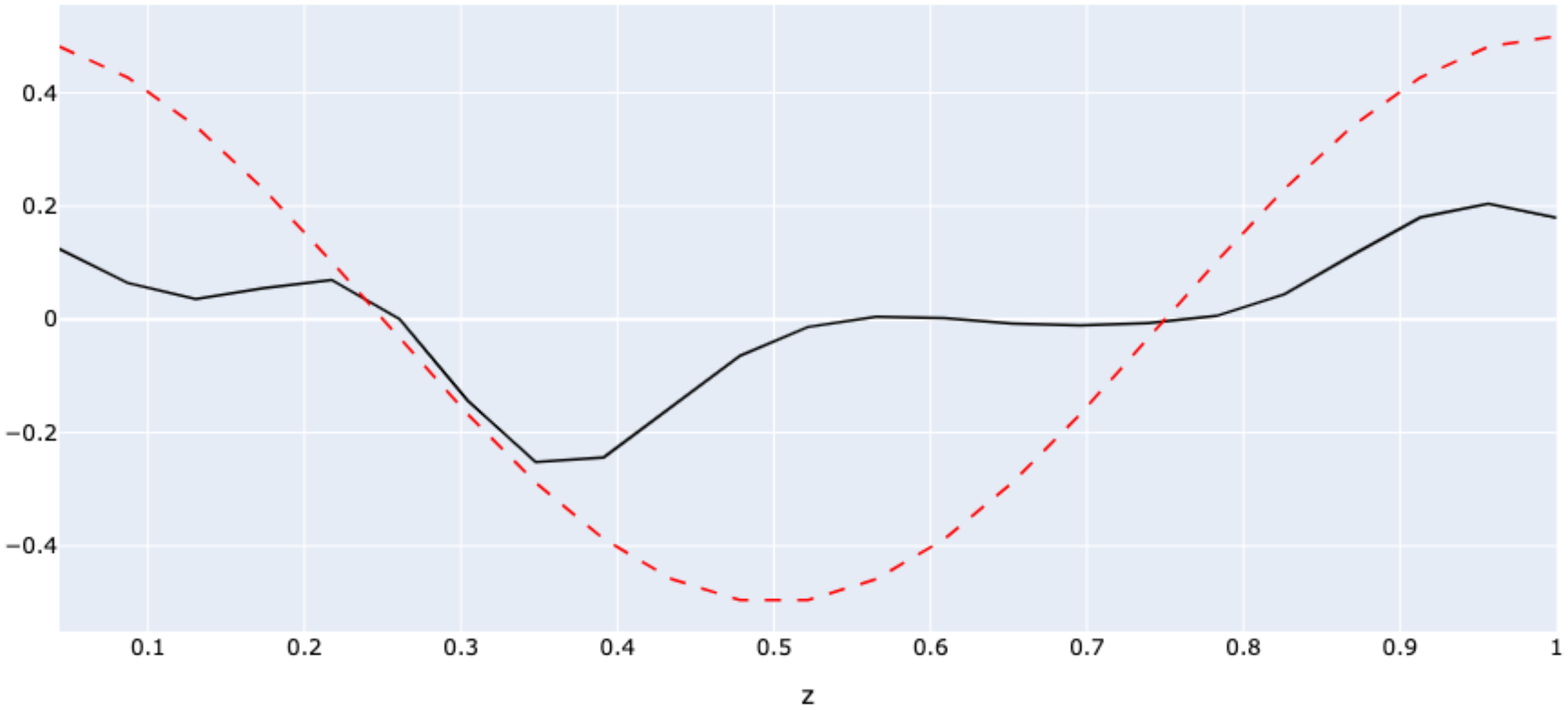}
\caption{\small Estimator of $S_\zs{1}$ for $n=500$}%
\label{fig.2-n-500}
\end{subfigure}
\caption{\small Model selection for $n=500$}\label{fig.n-500} 
\end{figure}


\newpage

\begin{figure}[h]\label{fig.n-10000}
\begin{subfigure}[b]{0.5\textwidth}

  \includegraphics[width=\columnwidth,height=6cm]{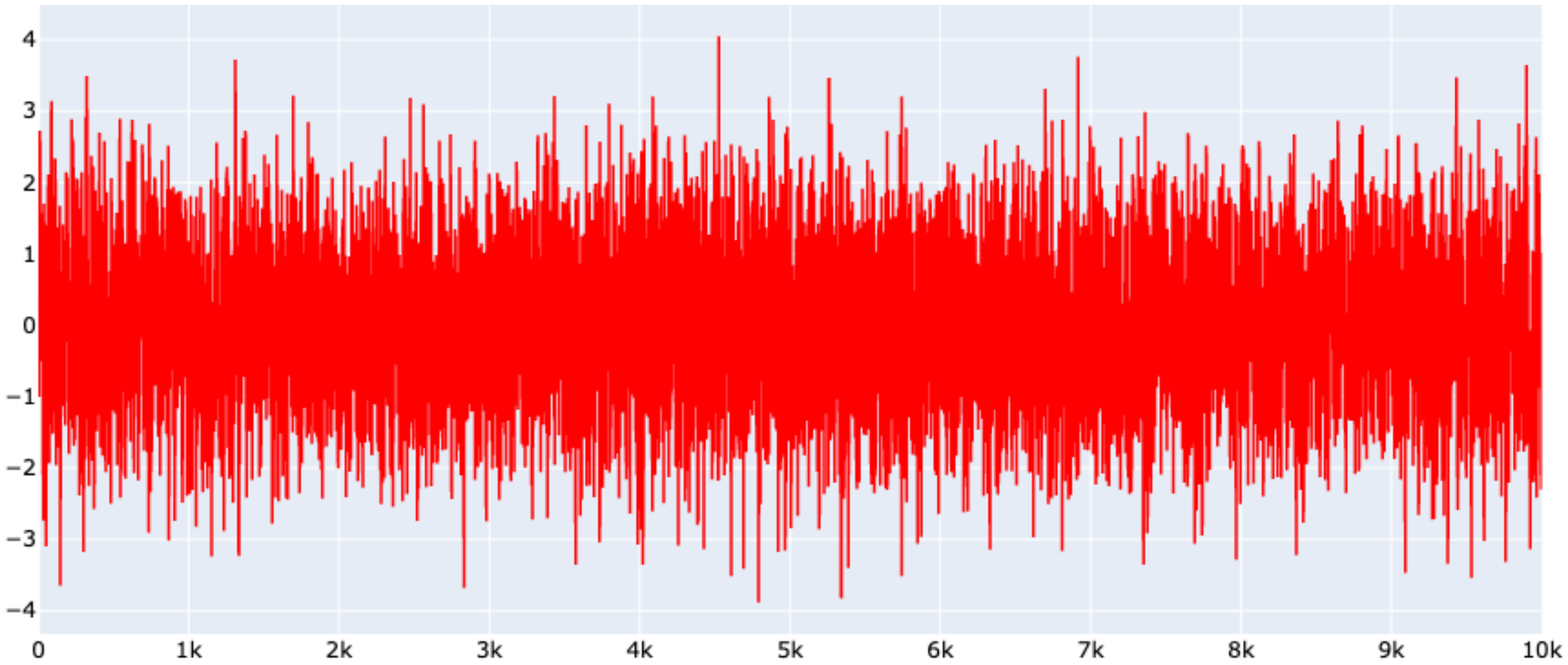}
\caption{\small Observations for $n=10000$}%
\label{fig.1-n-10000}
\end{subfigure}%
\begin{subfigure}[b]{0.5\textwidth}
\includegraphics[width=\columnwidth,height=6cm]{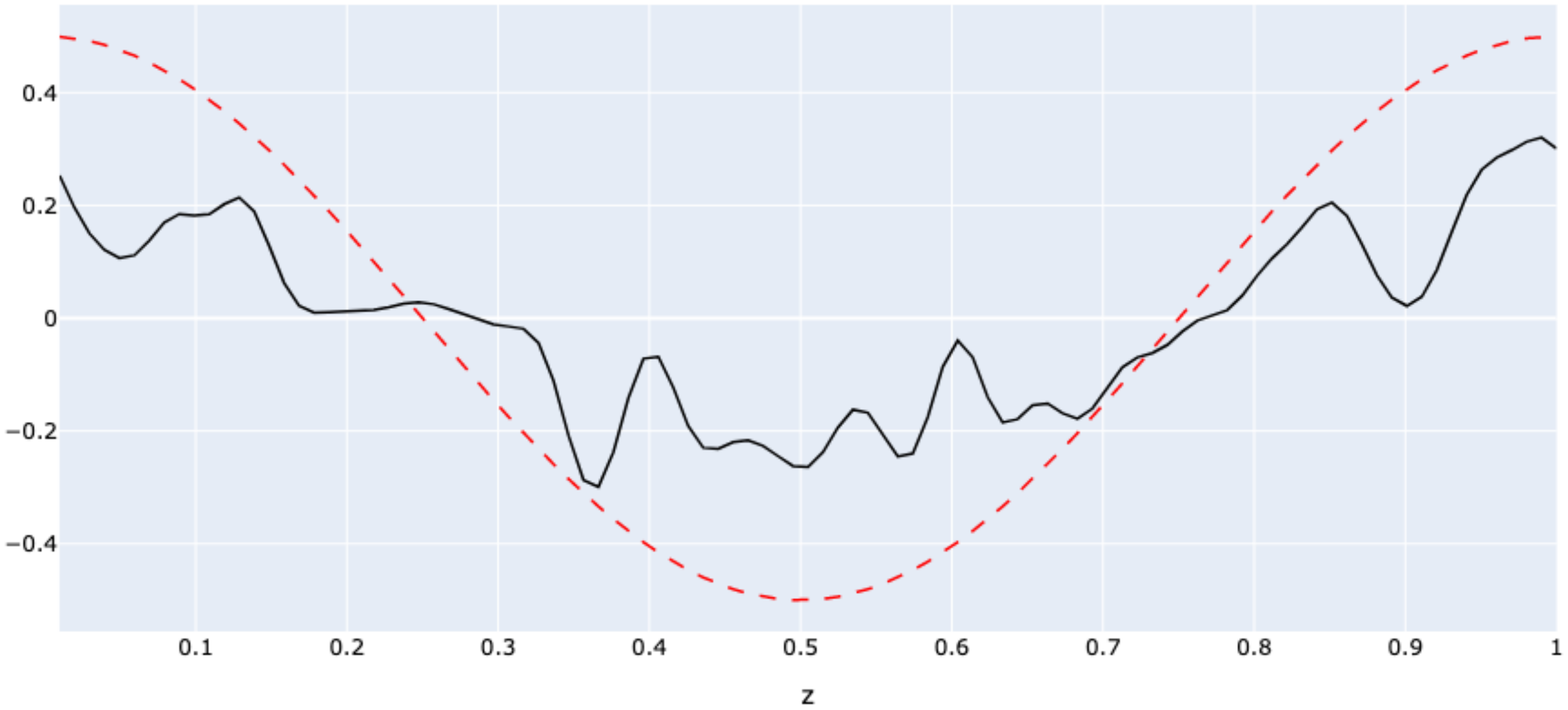}
\caption{\small Estimator of $S_\zs{1}$ for $n=10000$}%
\label{fig.2-n-10000}
\end{subfigure}
\caption{\small Model selection for $n=10000$}\label{fig.n-10000} 
\end{figure}


\begin{figure}[h]\label{fig.n-70000}
\begin{subfigure}[b]{0.5\textwidth}

 \includegraphics[width=\columnwidth,height=6cm]{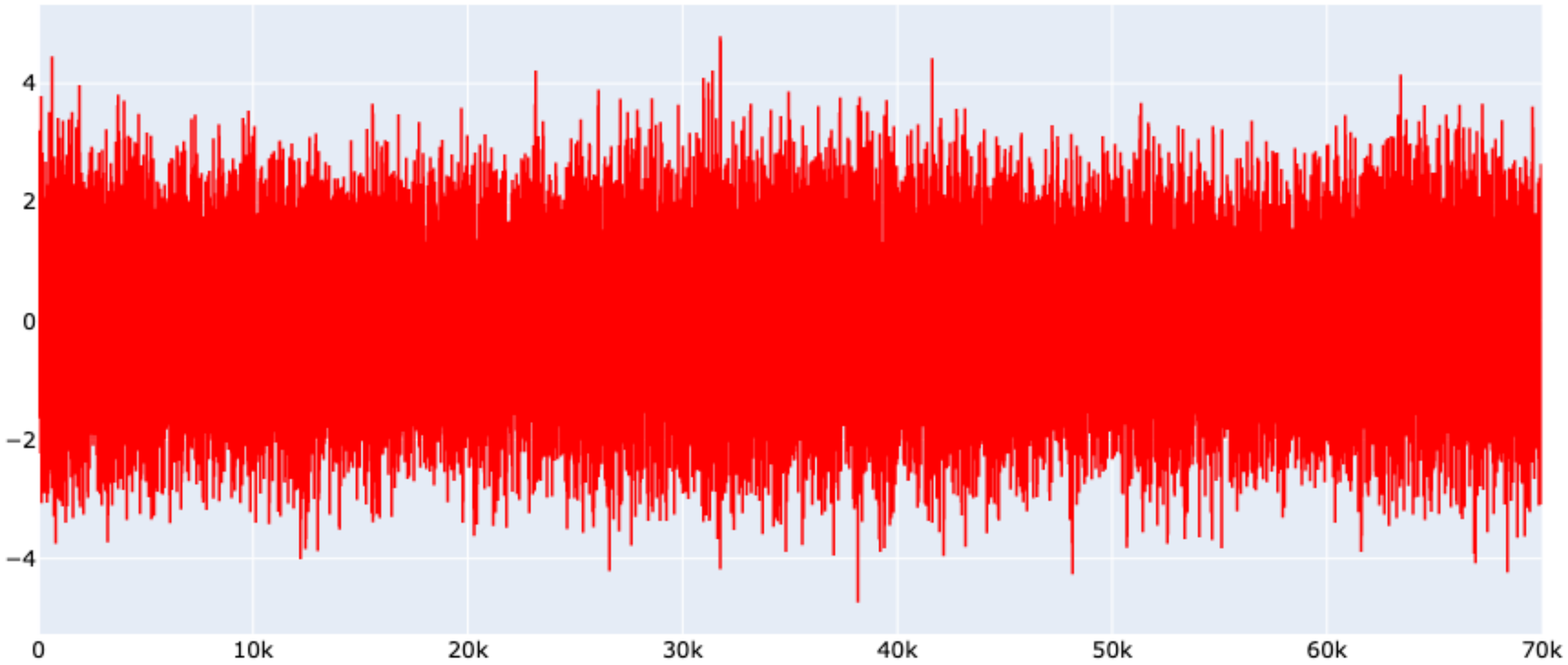}
\caption{\small Observations for $n=70000$}%
\label{fig.1-n-70000}
\end{subfigure}%
\begin{subfigure}[b]{0.5\textwidth}
\includegraphics[width=\columnwidth,height=6cm]{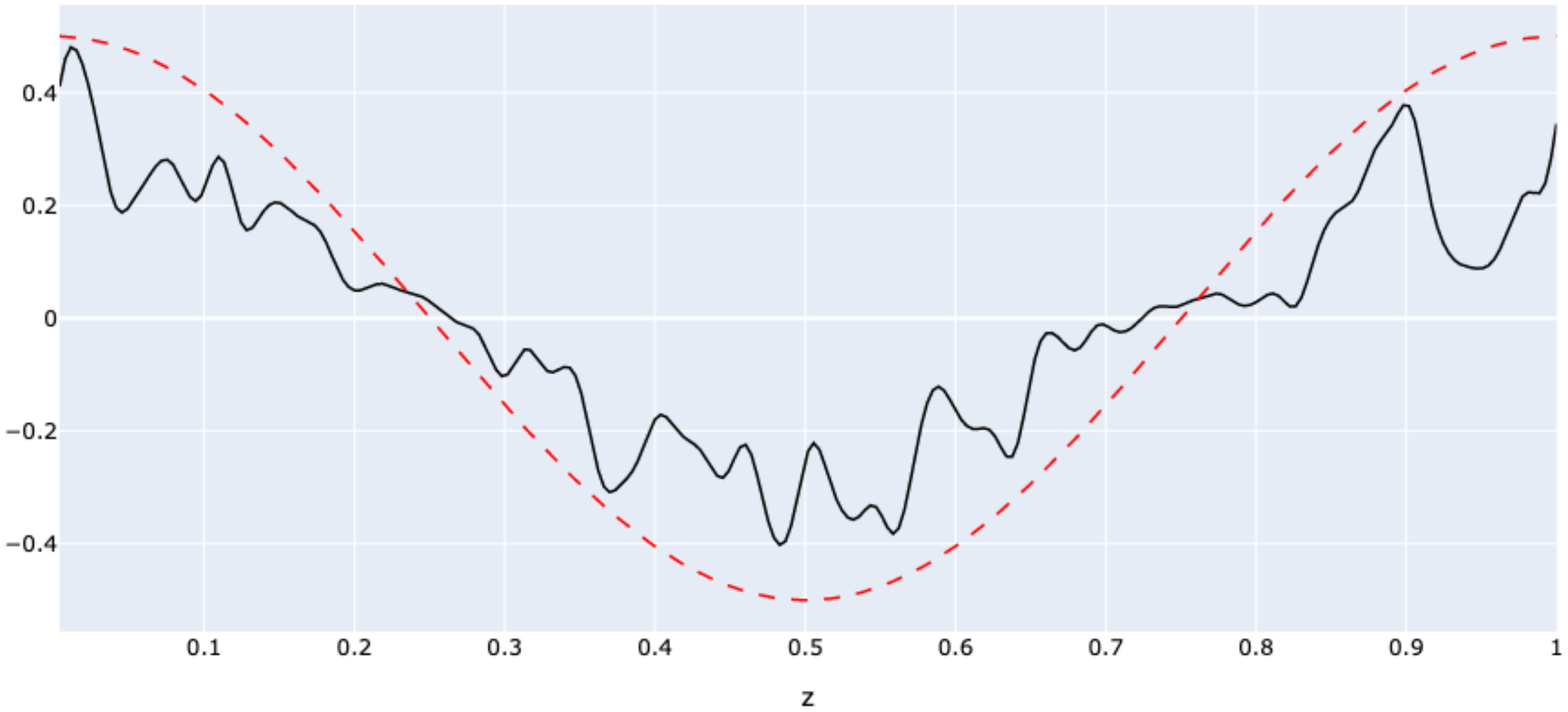}
\caption{\small Estimator of $S_\zs{1}$ for $n=70000$}%
\label{fig.2-n-70000}
\end{subfigure}
\caption{\small Model selection for $n=70000$}\label{fig.n-70000} 
\end{figure}


\newpage

Figures \ref{fig.n-200}--\ref{fig.n-70000} show  the behavior of the  function $S_\zs{1}$ and its estimators
by the model selection procedure \eqref{sec:Ms.13}
depending on the  observations number $n$.  In the figures (a) the observations are given and in (b)
the red dotted is the regression function and the black full line  is its  estimator at the points \eqref{sec:Dt.00-00}.
Then we calculate the empiric risks as
\begin{equation}\label{sec:Siml.1}
\overline{\cR}= \frac{1}{d} \sum_{j=1}^{d} \wh{\E} \left(\wh{S}_\zs{n}(z_\zs{j})-S(z_\zs{j})\right)^2
\,,
\end{equation}
where  the expectation is taken as an average over $M= 50$ replications, i.e.
$$
\wh{\E} \left(\wh{S}_\zs{n}(.)-S(.)\right)^2 = \frac{1}{M} \sum_\zs{l=1}^{M} \left(\wh{S}^l_\zs{n}(\cdot)-S(\cdot) \right)^2
\,.
$$
\noindent 
We use also the relative risk 
\begin{equation}\label{sec:Siml.2}
\overline{\cR}_\zs{*}=\frac{\overline{\cR}}{ \Vert  S\Vert  ^2_\zs{n}}
\quad\mbox{and}\quad
\Vert  S\Vert  ^2_n = \frac{1}{n} \sum_\zs{j=1}^n S^2(x_\zs{j})\,.
\end{equation}
The tables below give the values for the sample risks \eqref{sec:Siml.1} and \eqref{sec:Siml.2} for different numbers 
of observations $n$.

\begin{table}[h!]
\begin{center}
\caption{Empirical risks for $S_\zs{1}$} \label{tab:1}
    \begin{tabular}{|l|c|r|}
                                  \hline
                                  n & $\overline{\cR}$  & $\overline{\cR}_\zs{*}$  \\
                                  \hline
                                  200 &  0.135  & 0.98 \\
                                  \hline
                                  500 & 0.0893 & 0.624  \\
                                  \hline
                                  10000 & 0.043 & 0.362  \\
                                  \hline
                                  70000 & 0.03523  & 0.281  \\
                                  \hline
\end{tabular}
\end{center}
\end{table}

\begin{table}[h!]
\begin{center}
\caption{Empirical risks for $S_\zs{2}$} \label{tab:1}
    \begin{tabular}{|l|c|r|}
                                  \hline
                                  n & $\overline{\cR}$  & $\overline{\cR}_\zs{*}$  \\
                                  \hline
                                  200 &  0.0821  & 5.685 \\
                                  \hline
                                  500 & 0.0386 & 2.623  \\
                                  \hline
                                  10000 & 0.0071 & 0.516  \\
                                  \hline
                                  70000 & 0.0067  & 0.419  \\
                                  \hline
\end{tabular}
\end{center}
\end{table}

\begin{remark}\label{Re.sec:Mnc.1}
From numerical simulations  of the procedure \eqref{sec:Ms.13}
with various observation durations $n$ and for the different functions $S$ we may conclude  that the quality of the estimation improves as the number of observations increases.
\end{remark}

\newpage

\section{The van Trees inequality}\label{sec:VanTrees}

In this section we consider the nonparametric autoregressive model \eqref{sec:In.1}
with the  $(0,1)$ gaussian i.i.d. random variable
$(\xi_\zs{l})_\zs{1\le l\le n}$ and the parametric linear
 function $S$, i.e.
\begin{equation}\label{sec:VanTrees-1-S}
S_\zs{\theta}(x)=
\sum^{d}_\zs{j=1}\,\theta_\zs{j}\,\psi_\zs{j}(x)
\,,\quad \theta=(\theta_\zs{1},\ldots,\theta_\zs{d})'\in\bbr^{d}\,.
\end{equation}
\noindent
We assume that the functions $(\psi_\zs{j})_\zs{1\le j\le d}$ are orthogonal with respect to the scalar product
\eqref{sec:Ms.4nn_1}. Let now $\P^{n}_\zs{\theta}$ be the distribution in $\bbr^{n}$ of the  observations $y=(y_\zs{1},\ldots,y_\zs{n})$
in the model \eqref{sec:In.1}
with the function \eqref{sec:VanTrees-1-S}
and $\nu^{n}_\zs{\xi}$ be the distribution 
in
$\bbr^{n}$
of the gaussian vector
$(\xi_\zs{1},\dots,\xi_\zs{n})$.
In this case the Radon - Nykodim density  is given as 
\begin{equation}\label{sec:App.7}
f_\zs{n}(y,\theta)=
\frac{\d\P^{(n)}_\zs{\theta}}{\d\nu^{n}_\zs{\xi}}
=
\exp\left\{
\sum^{n}_\zs{l=1}S_\zs{\theta}(x_\zs{l})y_\zs{l-1}y_\zs{l}
-\frac{1}{2}
\sum^{n}_\zs{j=1}S^{2}_\zs{\theta}(x_\zs{l})y^{2}_\zs{l-1}
\right\}
\,.
\end{equation}

\noindent
Let $\varrho$ be a prior distribution density  on $\bbr^d$ for the parameter $\theta$ of
the following form
$
\varrho(\theta)=\prod_{j=1}^{d} \varrho_\zs{j}(\theta_\zs{j})$,
where $\varrho_\zs{j}$ is a some  probability density in $\bbr$ with continuously derivative $\dot{\varrho}_\zs{j}$ for which the Fisher information is finite, i.e.
\begin{equation}
\label{sec:FisherInf-1}
I_\zs{j}=\int_\zs{\bbr}\,\frac{\dot{\varrho}^2_\zs{j}(z)}{\varrho_\zs{j}(z)}\,\d z
<\infty
\,.
\end{equation}
\noindent 
Let $g(\theta)$ be a continuously differentiable $\bbr^{d}\to \bbr$ function such that
\begin{equation}\label{sec:App.8}
\lim_\zs{|\theta_\zs{j}|\to \infty}\,
g(\theta)\,\varrho_\zs{j}(\theta_\zs{j})=0
\quad\mbox{and}\quad
\int_\zs{\bbr^d}\,|g^{\prime}_\zs{j}(\theta)|\,\u(\theta)\,\d \theta
<\infty\,,
\end{equation}
where $g^{\prime}_\zs{j}(\theta)=\,\partial g(\theta)/\partial\theta_\zs{j}$.

\noindent
For any $\cB(\bbr^{n+d})$-measurable integrable function $H=H(y,\theta)$ we denote
\begin{equation}\label{Bayes-Rsk-1}
\wt{\E}\,H=
\int_\zs{\bbr^{n+d}}\,
H(y,\theta)\,f_\zs{n}(y,\theta)\,\varrho(\theta)
\d \nu^{(n)}_\zs{\xi}
\, \d \theta\,.
\end{equation}

\noindent 
Let  $\cF^y_\zs{n}$ be the field generated by the observations \eqref{sec:In.1}, i.e.
 $\cF^y_\zs{n}=\sigma\{y_\zs{1},\ldots,y_\zs{n}\}$.
Now we study the Gaussian the model \eqref{sec:In.1} with the function
\eqref{sec:VanTrees-1-S}.

\begin{lemma}\label{Le.sec:App.3}
For any $\cF^y_n$-measurable square integrable $\bbr^{n}\to\bbr$ function $\wh{g}_\zs{n}$
 and for any $1\le j\le d$,
 the mean square accuracy of the function $g(\cdot)$
 with respect to the distribution   \eqref{Bayes-Rsk-1}
 can be estimated from below as
\begin{equation}
\label{VanTreesInqly-1}
\wt{\E}(\wh{g}_\zs{n}-g(\theta))^2\ge
\frac{\overline{g}^{2}_\zs{j}}{
\wt{\E}\,
\Psi_\zs{n,j}
+I_\zs{j}}\,,
\end{equation}
where $\Psi_\zs{n,j}=\sum^{n}_\zs{l=1}\,
\psi^{2}_\zs{j}(x_\zs{l})
\,y^{2}_\zs{l-1}$ and
$\overline{g}_\zs{j}=\int_\zs{\bbr^{d}}\,g^{\prime}_\zs{j}(\theta)\,\varrho(\theta)\,\d \theta$.
\end{lemma}
\noindent {\bf  Proof.}
First, for any $\theta\in\bbr^{d}$ we set
$$
\wt{U}_\zs{j}=\wt{U}_\zs{j}(y,\theta)=
\frac{1}{f_\zs{n}(y,\theta)\u(\theta)}
\,
\frac{\partial\,(f_\zs{n}(y,\theta)\u(\theta))}{\partial\theta_\zs{j}}
 \,.
$$
Taking into account the condition \eqref{sec:App.8} and
integrating by parts we get
\begin{align*}
\wt{\E}&\left((\wh{g}_\zs{n}-g(\theta))\wt{U}_\zs{j}\right)
=\int_\zs{\bbr^{n+d}}\,(\wh{g}_\zs{n}(y)-g(\theta))\frac{\partial}{\partial\theta_\zs{j}}
\left(f_\zs{n}(y,\theta)\u(\theta)\right)\d \theta\,\d\nu^{(n)}_\zs{\xi}\\[2mm]
&=\int_\zs{\bbr^{n+d-1}}\left(\int^{+\infty}_\zs{-\infty}\,
g^{\prime}_\zs{j}(\theta)\,
f_\zs{n}(y,\theta)\u(\theta)\d \theta_\zs{j}\right)\left(\prod_\zs{i\neq j}\d \theta_i\right)\,\d \nu^{(n)}_\zs{\xi}
=\overline{g}_\zs{j}\,.
\end{align*}
\noindent
Now by the Bouniakovskii-Cauchy-Schwarz inequality
we obtain the following lower bound for the quadratic risk
$
\wt{\E}(\wh{g}_\zs{n}-g(\theta))^2\ge \overline{g}^{2}_\zs{j}/\wt{\E}\wt{U}_\zs{j}^2$.
To study the denominator in this
inequality note that in view of the representation
 \eqref{sec:App.7}
$$
\frac{1}{f_\zs{n}(y,\theta)}
\frac{\partial\,f_\zs{n}(y,\theta)}{\partial\theta_\zs{j}}
=
\sum^{n}_\zs{l=1}\,
\psi_\zs{j}(x_\zs{l})
\,y_\zs{l-1}
(
y_\zs{l}
-
S_\zs{\theta}(x_\zs{l})
y_\zs{l-1}
)
\,.
$$
Therefore, for each $\theta\in\bbr^{d}$,
$$
\E^{(n)}_\zs{\theta}\,
\frac{1}{f_\zs{n}(y,\theta)}
\frac{\partial\,f_\zs{n}(y,\theta)}{\partial\theta_\zs{j}}
\,
=0
$$
and
$$
\E^{(n)}_\zs{\theta}\,
\left(
\frac{1}{f_\zs{n}(y,\theta)}
\frac{\partial\,f_\zs{n}(y,\theta)}{\partial\theta_\zs{j}}
\right)^2
=\,
\E^{(n)}_\zs{\theta}
\sum^{n}_\zs{l=1}
\,
\psi^{2}_\zs{j}(x_\zs{l})
\,y^{2}_\zs{l-1}
=\,
\E^{(n)}_\zs{\theta}\,
\Psi_\zs{n,j}
\,.
$$
Using the equality
$$
\wt{U}_\zs{j}=
\frac{1}{f_\zs{n}(y,\theta)}
\frac{\partial\,f_\zs{n}(y,\theta)}{\partial\theta_\zs{j}}
+
\frac{1}{\u(\theta)}
\frac{\partial\,\u(\theta)}{\partial\theta_\zs{j}}
 \,,
$$
we get
$
\wt{\E}\,\wt{U}_\zs{j}^2=\,
\wt{\E}\,
\Psi_\zs{n,j}
+\,I_\zs{j}$.
Hence
Lemma~\ref{Le.sec:App.3}.
\endproof

\section{Lower bound}\label{sec:LowBnds-121}

First, taking into account that the $(0,1)$ gaussian density $\p_\zs{0}$ 
belongs to the class \eqref{2.1}, we get
 $
\cR^{*}(\wh{S}_\zs{n},S)\,
\ge\,
\cR_\zs{\p_\zs{0}}(\wh{S}_\zs{n},S)$.
Now, according to the general lower bounds methods  (see, for example, in \cite{GaltchoukPergamenshchikov2009b})
one needs to estimate this risk from below by some bayesian risk and, then to apply the van Trees inequlity.
To define the bayesian risk we need to choose a prior distribution on $\cW_\zs{\k,\r}$. To this end, first 
for any vector $\kappa=(\kappa_\zs{j})_\zs{1\le j\le d}\in\bbr^{d}$,
 we set
\begin{equation}\label{sec:Lo.4--1}
S_\zs{\kappa}(x)=\sum_{j=1}^{d_\zs{n}}\,\kappa_\zs{j}\,\phi_\zs{j}(x)
\,,
\end{equation}
where $(\phi_\zs{j})_\zs{1\le j\le d_\zs{n}}$ is the trigonometric basis defined in \eqref{sec:In.5_trg}.
 We will choose a prior distribution on the basis of 
 the optimal coefficients in  \eqref{sec:Lo.4--1} (see, for example, in \cite{PchelintsevPergamenshchikovPovzun2021, Pinsker1981})
 providing asymptotically the maximal valuer for the risk $\cR_\zs{\p_\zs{0}}(\wh{S}_\zs{n},S)$ over $S$.  So, we 
 choose it as the distribution $\mu_\zs{\kappa}$ in $\bbr^{d}$ of the random vector
$\kappa=(\kappa_\zs{j})_\zs{1\le j\le d_\zs{n}}$ defined by its components as
\begin{equation}\label{sec:Lo.5}
\kappa_\zs{j}=s_\zs{j}\,\eta^{*}_\zs{j}
\,,
\end{equation}
where $\eta^{*}_\zs{j}$ are i.i.d. random variables with the continuously  differentiable density $\rho_\zs{n}(\cdot)$ defined in
Lemma \ref{Le.sec:App.PriDens-1} with $\N=\ln n$ and $n> e^{2}$, 
$$
s_\zs{j}=\sqrt{\frac{(b-a)s^*_\zs{j}}{n}}
\quad\mbox{and}\quad
s^{*}_\zs{j}\,
=
\left( \frac{d_\zs{n}}{j}
\right)^{\k}
-
1
\,.
$$
\noindent
To choose the number of the terms $d$ in \eqref{sec:Lo.4--1} one needs to keep this function in $\cW_\zs{\k,\r}$, i.e. one needs to provide
for arbitrary fixed $0<\rho<1$
 the following property
\begin{equation}\label{sec:Lo.PR-d}
\lim_\zs{n\to\infty}\sum^{d}_\zs{j=1}\a_\zs{j}\,s^{2}_\zs{j}=\rho \r:=\r_\zs{\rho}\,.
\end{equation}

\noindent
To do this we set
\begin{equation}\label{sec:Lo.1}
d=d_\zs{n}=\left[\g^{*}_\zs{\k}\, n^{1/(2\k+1)}\right]\,,
\end{equation}
where $\g^{*}_\zs{\k}=(b-a)^{(2\k-1)/(2\k+1)}\l_\zs{\k}(\r_\zs{\rho})(\k+1)/\k$
 and
$$
\l_\zs{\k}(\r_\zs{\rho})
=\,((1+2\k)\r_\zs{\rho})^{1/(2\k+1)}\,
\left(\frac{\k}{\pi (\k+1)}\right)^{2\k/(2\k+1)}
=\rho^{1/(2\k+1)}\l_\zs{\k}(\r)
\,.
$$

\noindent
It is clear that almost sure the function \eqref{sec:Lo.4--1} can be bounded as
\begin{equation}
\label{UnifUpperBound-S}
\max_\zs{a\le x\le b}
\left(
\vert S_\zs{\kappa}(x)\vert
+
\vert \dot{S}_\zs{\kappa}(x)\vert
\right)
\le 
\c_\zs{*}
\frac{\ln n}{\sqrt{n}}
\sum^{d_\zs{n}}_\zs{j=1}j\left(\frac{d_\zs{n}}{j}\right)^{\k/2}
:=\delta^{*}_\zs{n}\,,
\end{equation}
where
$
\c_\zs{*}
=\sqrt{2}(b-a+\pi)/(b-a)^{3/2}$.
Note, that  $\delta^{*}_\zs{n}\to 0$ as $n\to\infty$ for  $\k\ge 2$. 
Therefore, for sufficiently large $n$  the function  \eqref{sec:Lo.4--1}
belongs to  the class \eqref{sec:Sp.1}. Now, $\forall f\in\cL_\zs{2}[a,b]$, we denote by $\h(f)$ its projection
in $\cL_\zs{2}[a,b]$
 onto the ball
 $\cW^{*}_\zs{\r}=\{f\in\cL_\zs{2}[a,b]\,:\,\Vert f\Vert\le \r\}$, i.e.
 $\h(f)=(\r/\max(\r\,,\,\Vert f\Vert)) f$.
Since  $\cW_\zs{\k,\r}\subset \cW^{*}_\zs{\r}$, then for $S\in \cW_\zs{\k,\r}$ we have
$
\|\wh{S}_\zs{n}-S\|^2\ge\|\wh{\h}_\zs{n}-S\|^2$,
where
$\wh{\h}_\zs{n}=\h(\wh{S}_\zs{n})$.
Therefore,  for large $n$
\begin{align}\nonumber
\sup_\zs{S\in \cW_\zs{\k,\r}}\,\upsilon(S)\,\cR_\zs{\p_\zs{0}}(\wh{S}_\zs{n},S)\,
&\ge\,
\int_\zs{\D_\zs{n}}\,
\upsilon(S_\zs{z})
\E_\zs{\p_\zs{0},S_\zs{z}}\|\wh{\h}_\zs{n}-S_\zs{z}\|^2\,\mu_{\kappa}(\d z)\\[2mm]\label{LB-1-BR-01}
&\ge\,
\upsilon_\zs{*}\,
\int_\zs{\D_\zs{n}}\,
\E_\zs{\p_\zs{0},S_\zs{z}}\|\wh{\h}_\zs{n}-S_\zs{z}\|^2\,\mu_{\kappa}(\d z)
\,,
\end{align}
where 
$
\D_\zs{n}=\left\{z\in\bbr^{d}\,:\,\sum^{d}_\zs{j=1}\,\a_\zs{j} z^{2}_\zs{j}\le \r
\right\}$
and
$
\upsilon_\zs{*}
=\inf_\zs{\vert S\vert_\zs{*} \le \delta^{*}_\zs{n}}\,\upsilon(S)$.
Note that
$\upsilon_\zs{*}
\to (b-a)^{-4\k/(2\k+1)}$ as
$
 n\to\infty$.
Using the distribution $\mu_\zs{\kappa}$ we introduce
 the following Bayes risk as
\begin{equation}
\label{sec:Bayes-risk-1}
\wt{\cR}_\zs{0}(\wh{S}_\zs{n})=
\wt{\E}_\zs{0}\,\Vert\wh{S}_\zs{n}-S_\zs{z}\Vert^{2}
=
\int_\zs{\bbr^d}
\E_\zs{\p_\zs{0},S_\zs{z}}
\,
\Vert\wh{S}_\zs{n}-S_\zs{z}\Vert^{2}
\,
\mu_\zs{\kappa}(\d z)
\,.
\end{equation}
\noindent
Now taking into account that $\|\wh{\h}_\zs{n}\|^2\le r$,
 we get
\begin{equation}\label{sec:Lo.12}
\sup_\zs{S\in \cW_\zs{\k,\r}}\,
\upsilon(S)
\,
\cR_\zs{\p_\zs{0}}(\wh{S}_\zs{n},S)
\,
\ge\,\upsilon_\zs{*}\,
\wt{\cR}_\zs{0}(\wh{\h}_\zs{n})
-2\,\upsilon_\zs{*}\,
\R_\zs{0,n}
\end{equation}
and
$
\R_\zs{0,n}=
\int_\zs{\D^{c}_\zs{n}}\,
\,
(r+\|S_\zs{z}\|^2)\,
\mu_\zs{\kappa}(\d z)
=
\int_\zs{\D^{c}_\zs{n}}\,
\left(\r+\,\vert z\vert^2\right)\,
\mu_\zs{\kappa}(\d z)$.
Note here that this term is studied In
Lemma \ref{Le.sec:App.3+1}. Moreover,
note also that  for any $z\in\bbr^{d}$
we get
$
\|\wh{\h}_\zs{n}-S_z\|^2 \ge
\sum_{j=1}^{d_\zs{n}}\,
(\wh{z}_\zs{j}-z_\zs{j})^2$
and
$\wh{z}_\zs{j}=\left(\wh{\h}_\zs{n}\,,\phi_\zs{j}\right)$. Therefore, from  Lemma~\ref{Le.sec:App.3} 
 with $g(\theta)=\theta_\zs{j}$  it follows, that for any $1\le j\le d$
and any $\cF^{y}_\zs{n}$ measurable random variable $\wh{\kappa}_\zs{j}$
$$
\wt{\E}_\zs{0}(\wh{\kappa}_\zs{j}-\kappa_\zs{j})^{2}
\ge 
\frac{1}{\wt{\E}_\zs{0}\,\Psi_\zs{n,j}
+s^{-2}_\zs{j} J_\zs{n}}\,,
 $$
 where
 $
\Psi_\zs{n,j}=\sum^{n}_\zs{l=1}\phi^{2}_\zs{j}(x_\zs{l})y^{2}_\zs{l-1} $
and
$
J_\zs{n}=\int^{\ln n}_\zs{-{\ln n}}\,(\dot{\rho}_\zs{n}(t))^{2}/\rho_\zs{n}(t)
\d t$.
Therefore, the  Bayes risk can be estimated from below as
 $$
\wt{\cR}_\zs{0}(\wh{\h}_\zs{n})\,
\ge\,
\sum_{j=1}^{d_\zs{n}}\,\frac{1}{\wt{\E}_\zs{0}\,
\Psi_\zs{n,j}+s^{-2}_\zs{j}J_\zs{n}}
=\frac{b-a}{n} 
\sum_{j=1}^{d_\zs{n}}\,\frac{1}{
\overline{\Psi}_\zs{n,j}+(s^{*}_\zs{j})^{-1}J_\zs{n}}
 \,,
 $$
where 
\begin{equation}
\label{Psi-Def-1}
\overline{\Psi}_\zs{n,j}=
\frac{(b-a)}{n}
\wt{\E}_\zs{0}\,\Psi_\zs{n,j}
=
\frac{(b-a)}{n}
\sum^{n}_\zs{l=1}\phi^{2}_\zs{j}(x_\zs{l})\wt{\E}_\zs{0}\,y^{2}_\zs{l-1} 
\,.
\end{equation}
\noindent
Note here, that in view of  Lemmas \ref{Le.sec:App.PriDens-1} and \ref{Le.sec:App.PriDstr-00}
  for any $0<\rho_\zs{1}<1$ and sufficiently large  $n$ we have $J_\zs{n}\le 1+\rho_\zs{1}$ and
  $ \max_\zs{1\le j\le d}\overline{\Psi}_\zs{n,j}\le 1+\rho_\zs{1}$, therefore,
$$
\wt{\cR}_\zs{0}(\wh{\h})\,
\ge\,
\frac{b-a}{(1+\rho_\zs{1})n}\,
\sum_{j=1}^{d_\zs{n}}\,\frac{s^{*}_\zs{j}}{1+s^{*}_\zs{j}}
=
\,\frac{b-a}{(1+\rho_\zs{1})n}\,
\sum_\zs{j=1}^{d_\zs{n}}\,
\left(
1
-
\frac{j^\k}{d^\k_\zs{n}}
\right)
\,.
$$
\noindent
Using here that
$
\lim_\zs{d\to\infty} d^{-1}
\sum_\zs{j=1}^{d}\,
\left(
1
- j^{\k} d^{-\k}
\right)
=\int^{1}_\zs{0}(1-t^{\k})\d t=\k/(\k+1)
$, 
we obtain for sufficiently large $n$
$$
\wt{\cR}_\zs{0}(\wh{\h})\,
\ge
\,\frac{(b-a)(1-\rho_\zs{1})d_\zs{n}}{(1+\rho_\zs{1})n}\,\frac{\k}{\k+1}
\,.
$$
\noindent
Therefore, taking into account this  in
\eqref{sec:Lo.12}
we conclude
  through the definition \eqref{sec:Lo.1}
and
Lemma \ref{Le.sec:App.3+1}, that for any   $0<\rho$ and $\rho_\zs{1}<1$
$$
\liminf_\zs{n\to\infty}\inf_\zs{\wh{S}_\zs{n}\in\Xi_\zs{n}}\,n^{\frac{2k}{2k+1}}\,
\sup_\zs{S\in \cW_\zs{\k,\r}}\,\upsilon(S)\,\cR^{*}(\wh{S}_\zs{n},S)
\ge\,
\frac{\rho^{\frac{1}{2k+1}}(1-\rho_\zs{1})}{1+\rho_\zs{1}}\,
\l_\zs{\k}(\r)\,.
$$
Taking here limit as $\rho\to 1$ and  $\rho_\zs{1}\to 0$ we come to the Theorem~\ref{Th.2.1-1}.
\endproof

\section{Upper bound}\label{sec:MainResUp}

 We start with the estimation problem 
for the functions $S$ from $\cW_\zs{\k,\r}$ with known
 parameters $k$, $r$ and $\varsigma_\zs{*}$ defined in \eqref{2.7}.
In this case we use the estimator from family \eqref{sec:Ga.3}
\begin{equation}\label{U.1}
\wt{S}=\wh{S}_\zs{\wt{\alpha}}
\quad\mbox{and}\quad
\wt{\alpha}=(\k,\wt{\t}_n)
\,,
\end{equation}
 where  
$\wt{\t}_\zs{n}=
\left[\overline{r}(S)/\varepsilon
\right] \varepsilon$,
$\overline{r}(S)=\r/\varsigma_\zs{*}$ and   $\varepsilon=1/\ln n$.
Note that for sufficiently large $n$,
the parameter 
$\wt{\alpha}$ belongs to the set \eqref{sec:Ga.0}.
In this section we study the risk for   the estimator \eqref{U.1}. To this end we need firstly to analyse the
 asymptotic behavior  of the sequance
 \begin{equation}
 \label{sec:U-UPsil-1}
\Upsilon_\zs{n}(S)=\sum_{j=1}^{d}(1-\wt{\lambda}(j))^2 \theta^2_\zs{j,d}+ \frac{\varsigma_\zs{*}}{n}
\,
\sum_{j=1}^{d}\,\wt{\lambda}^2(j)\,.
\end{equation}

\begin{proposition}\label{Pr.sec:App.44}
The sequence $\Upsilon_\zs{n}(S)$  is bounded from above
$$
\limsup_\zs{n\to\infty}\sup_\zs{S\in \cW_\zs{\k,\r}}\,n^{\k_\zs{1}}\,
\upsilon(S)\,\Upsilon_\zs{n}(S)\,\le\,\l_\zs{\k}(\r)
\quad\mbox{and}\quad
\k_\zs{1}=2\k/(2\k+1)
\,.
$$
\end{proposition}

\proof 
\noindent
First,  note that
$
0<\epsilon^{2} (b-a)\le\inf_\zs{S\in\Theta_\zs{\varepsilon,L}}\,\varsigma_\zs{*}
\le
\sup_\zs{S\in\Theta_\zs{\varepsilon,L}}\,\varsigma_\zs{*}
\le b-a$.
This implies directly that
\begin{equation}\label{U.9}
\lim_\zs{n\to\infty} \sup_\zs{S\in\Theta_\zs{\varepsilon,L}}\,
\left|
\wt{\t}_\zs{n}/\overline{r}(S)-1\right|=0\,,
\end{equation}
where $\wt{\t}_\zs{n}=
\left[\overline{r}(S)/\varepsilon
\right] \varepsilon$ and $\overline{r}(S)=\r/\varsigma_\zs{*}$.
\noindent Moreover, note that
$$
n^{\k_\zs{1}}
\upsilon(S)\,\Upsilon_\zs{n}(S)\le n^{\k_\zs{1}}\,\upsilon(S)\,\G_\zs{n}
+\,\frac{(\varsigma_\zs{*})^{1/(2\k+1)}}{(b-a)^{\k_\zs{1}} n^{1/(2\k+1)}}\,
\sum_{j=1}^{d}\,\wt{\lambda}^2(j)
$$
and
$
 \G_\zs{n}=
\sum^{d}_{j=1}\,(1-\wt{\lambda}(j))^2\,\theta^2_\zs{j,d}
=
\G_\zs{1,n}+\G_\zs{2,n}$,
where
$$
\G_\zs{1,n}=\sum_\zs{j=j_\zs{*}}^{[\wt{\omega}]}\,(1-\wt{\lambda}(j))^2\,\theta^2_\zs{j,d}
\quad\mbox{and}\quad
\G_\zs{2,n}=\sum_{j=[\wt{\omega}]+1}^{d}\,\theta^2_\zs{j,d}\,.
$$
Remind, that  
$
\wt{\omega}
=
\wt{\omega}_\zs{*}
+
(b-a)^{\k_\zs{1}}
 \left(
\varpi_\zs{\k}
 \wt{\t}_\zs{n}\,n
 \right)^{1/(2\k+1)}$
 and
$\varpi_\zs{\k}= (\k+1)(2\k+1)/(\pi^{2\k}\k)$.
\noindent
Note now, that
 Lemma~\ref{Le.sec:A.3-00}
 and Lemma~\ref{Le.sec:A.3}
  yield
$$
\G_\zs{1,n}
\le (1+\wt{\varepsilon})\,
\sum_\zs{j=j_\zs{*}}^{[\wt{\omega}]}\,(1-\wt{\lambda}(j))^2\,\theta^{2}_\zs{j}
+4
\r
(1+\wt{\varepsilon}^{-1})\,
\frac{(b-a)^{2\k}\wt{\omega}}{d^{2k}}
$$
and
$$
\G_\zs{2,n}\le (1+\wt{\varepsilon})
\sum_\zs{j\ge [\wt{\omega}]+1}\,\theta^{2}_\zs{j}
+
\r
(1+\wt{\varepsilon}^{-1})
\frac{(b-a)^{2\k}}{d^{2} 
\,\wt{\omega}^{2(k-1)}}
\,,
$$
i.e.
$
 \G_\zs{n}\le
 (1+\wt{\varepsilon})
 \G^{*}_\zs{n}
 +
 4\r (b-a)^{2\k}
 (1+\wt{\varepsilon}^{-1})\,
 \wt{\G}_\zs{n}$,
where
$$
\G^{*}_\zs{n}=
\sum_\zs{j\ge 1}(1-\wt{\lambda}(j))^2\,\theta^{2}_\zs{j}
=
\sum_\zs{j\le \wt{\omega}}(1-\wt{\lambda}(j))^2\,\theta^{2}_\zs{j}
+
\sum_\zs{j> \wt{\omega}}\,\theta^{2}_\zs{j}
:=
\G^{*}_\zs{1,n}+
\G^{*}_\zs{2,n}
$$
and
$
\wt{\G}_\zs{n}
= \wt{\omega}\,d^{-2\k}
 +d^{-2} \wt{\omega}^{-2(\k-1)}$.
Note, that
$$
n^{\k_\zs{1}}\upsilon(S)
\G^{*}_\zs{1,n}=\frac{\upsilon(S)}{(b-a)^{2\k\k_\zs{1}}(\varpi_\zs{\k}\wt{\t}_\zs{n})^{\k_\zs{1}}}\,
\sum_\zs{j=j_\zs{*}}^{[\wt{\omega}]}\,j^{2\k}\,\theta^{2}_\zs{j}\\[2mm]
\le 
\frac{\u^{*}_\zs{n}}{(\varpi_\zs{\k}\varsigma_\zs{*} \wt{\t}_\zs{n})^{\k_\zs{1}}}
\sum_\zs{j= j_\zs{*}}^{[\wt{\omega}]}\,\a_\zs{j}\,\theta^{2}_\zs{j}
\,,
$$
where $\u^{*}_\zs{n}=\sup_\zs{j\ge j_\zs{*}}\,j^{2k}/((b-a)^{2\k}\a_\zs{j})$. It is clear that $\lim_\zs{n\to\infty}\,\u^{*}_\zs{n}=\pi^{-2k}$.
Therefore, from
\eqref{U.9} it follows that
\begin{align*}
\limsup_\zs{n\to\infty}
\sup_\zs{S\in\Theta_\zs{\varepsilon,L}}\,
\frac{n^{\k_\zs{1}}\upsilon(S)\,\G^{*}_\zs{1,n}}{\sum_\zs{j=j_\zs{*}}^{[\wt{\omega}]}\,\a_\zs{j}\,\theta^{2}_\zs{j}}
\le
 \pi^{-2\k}(\varpi_\zs{\k}\r)^{-\k_\zs{1}}
\,.
\end{align*}
\noindent
Next note, that
for any $0<\wt{\varepsilon}<1$
and for sufficiently large $n$ 
\begin{align*}
\G^{*}_\zs{2,n}
\le \frac{1}{\a_\zs{[\wt{\omega}]+1}}
\sum_\zs{j\ge [\wt{\omega}]+1}\,\a_\zs{j}\,\theta^{2}_\zs{j}
\le
\frac{(1+\wt{\varepsilon})((b-a)\varsigma_\zs{*})^{\k_\zs{1}}}{\pi^{2\k}\,(\varpi_\zs{\k} \r n)^{\k_\zs{1}}}
\sum_\zs{j\ge [\wt{\omega}]+1}\a_\zs{j}\,\theta^{2}_\zs{j}\,,
\end{align*}
i.e.
$$
\limsup_\zs{n\to\infty}
\sup_\zs{S\in\Theta_\zs{\varepsilon,L}}\,
\frac{n^{\k_\zs{1}}\upsilon(S)\,\sum_\zs{j\ge [\wt{\omega}]+1}\,\theta^{2}_\zs{j}}{\sum_\zs{j=j_\zs{*}}^{[\wt{\omega}]}\a_\zs{j}\,\theta^{2}_\zs{j}}
\le 
\pi^{-2\k}(\varpi_\zs{\k}\r)^{-\k_\zs{1}}
\,.
$$
\noindent
Next, since $\sum_{j=1}^{d}\,\wt{\lambda}^2(j)/\wt{\omega}\to \k\k_\zs{1}/(\k+1)$ as $n\to\infty$, 
 we get directly that
$$
\lim_\zs{T\to\infty}\sup_\zs{S\in\Theta_\zs{\varepsilon,L}}
\left|
\frac{(\varsigma_\zs{*})^{1/(2\k+1)}
\sum_{j=1}^{d}\,\wt{\lambda}^2(j)}{(b-a)^{\k_\zs{1}}n^{1/(2\k+1)}}
-
\frac{(\varpi_\zs{\k}\r)^{1/(2\k+1)}\k \k_\zs{1} }{(\k+1)}
\right|=0\,.
$$
\noindent 
Finally, taking into account, that
$
\lim_\zs{n\to\infty}\,
\sup_\zs{S\in \cW_\zs{\k,\r}}
n^{\k_\zs{1}}
\wt{\G}_\zs{n}
=0$,
we obtain  
Proposition \ref{Pr.sec:App.44}.
\endproof

\begin{theorem}\label{Th.Up.1} 
The estimator $\wt{S}$ constructed on the trigonometric basis
satisfies the following asymptotic upper bound
\begin{equation}\label{U.2}
\limsup_\zs{n\to\infty}\,n^{2\k/(2\k+1)}\,
\sup_\zs{S\in \cW_\zs{\k,\r}}
\upsilon(S)\,\E_\zs{\p,S}\|\wt{S}-S\|^2_\zs{d}\,\Chi_\zs{\Gamma}\,\le \l_\zs{\k}(\r)\,.
\end{equation}
\end{theorem}
\noindent {\bf Proof.}  We denote $\wt{\lambda}=\lambda_\zs{\wt{\alpha}}$
and $\wt{\omega}=\omega_\zs{\wt{\alpha}}$.
Now we recall  that the Fourier coefficients 
on the set $\Gamma$
$$
\wh{\theta}_\zs{j,d}\,=\,\theta_\zs{j,d}\,+\,\zeta_\zs{j,d}
\quad\mbox{with}\quad
\zeta_\zs{j,d}=\sqrt{\frac{b-a}{d}}
\eta_\zs{j,d}+\varpi_\zs{j,d}\,.
$$
Therefore, on the set $\Gamma$
we can represent the empiric squared error as 
\begin{align*}
\|\wt{S}-S\|_\zs{d}^2
&=\sum_{j=1}^{d}\,(1-\wt{\lambda}(j))^2\,
\theta^2_\zs{j,d}-2 M_\zs{n}\\ 
&-
2\sum_{j=1}^{d}\,(1\,-\,\wt{\lambda}(j))\,\wt{\lambda}(j)
\theta_\zs{j,d}\,\varpi_\zs{j,d}
+\,
\sum_{j=1}^{d}\,\wt{\lambda}^2(j)\,\zeta^2_\zs{j,d}\,,
\end{align*}
where
$
M_\zs{n}
=
\sqrt{b-a}
\sum_{j=1}^{d}\,(1\,-\,\wt{\lambda}(j))\,\wt{\lambda}(j)
\theta_\zs{j,d}\,\eta_\zs{j,d}/\sqrt{d}$.
Now for any $0<\varepsilon_\zs{1}<1$ 
\begin{align*}
2\left|
\sum_{j=1}^{d}(1-\wt{\lambda}(j))\wt{\lambda}(j) \theta_\zs{j,d}\,\varpi_\zs{j,d}
\right|\,
&\le\, 
\varepsilon_\zs{1}\,\sum_{j=1}^{d}(1-\wt{\lambda}(j))^2\theta^2_\zs{j,d}\,
+
\varepsilon_\zs{1}^{-1}\sum_{j=1}^{d}\varpi^2_\zs{j,d}\\
&\le 
x
\varepsilon_\zs{1}\,\sum_{j=1}^{d}(1-\wt{\lambda}(j))^2\theta^2_\zs{j,d}\,
+
\frac{\varpi^{*}_\zs{n}}{\varepsilon_\zs{1} n}
\,,
\end{align*}
where $\varpi^{*}_\zs{n}$ is defined in \eqref{sec:Ms.6}.
Therefore,
$$
\|\wt{S}-S\|_\zs{d}^2\,
\le\, (1+\varepsilon_\zs{1})\sum_{j=1}^{d}
(1-\wt{\lambda}(j))^2 \theta^2_\zs{j,d}
-2M_\zs{n}
+\frac{\varpi^{*}_\zs{n}}{\varepsilon_\zs{1} n}
+
\sum_{j=1}^{d}
\wt{\lambda}^2(j)
\zeta^2_\zs{j,d}\,.
$$
By the same way we get
$$
\sum_{j=1}^{d}\,\wt{\lambda}^2(j)\,\zeta^2_\zs{j,d}
\le\,
\frac{(1+\varepsilon_\zs{1})(b-a)}{d}
\sum_{j=1}^{d}\,\wt{\lambda}^2(j)\,\eta^2_\zs{j,d}
+(1+\varepsilon_\zs{1}^{-1})\frac{\varpi^{*}_\zs{n}}{n}\,.
$$
Thus, on the set $\Gamma$ we find that for any $0<\varepsilon_\zs{1}<1$
$$
\|\wt{S}_\zs{n}-S\|_\zs{d}^2
\le
(1+\varepsilon_\zs{1})\Upsilon_\zs{n}(S)
-2M_\zs{n}
+(1+\varepsilon_\zs{1})\U_\zs{n}
+\frac{3\varpi^{*}_\zs{n}}{\varepsilon_\zs{1} n}\,,
$$
where $\Upsilon_\zs{n}(S)$ is defined in \eqref{sec:U-UPsil-1}
and
\begin{equation}\label{U.7-010}
\U_\zs{n}=
\frac{1}{d^{2}}\sum_{j=1}^{d}
\wt{\lambda}^2(j)
\left(d (b-a)
\eta^2_\zs{j,d}-\varsigma_\zs{*}
\right)
\,.
\end{equation}
We recall that the variance $\varsigma_\zs{*}$ is defined in \eqref{2.7}.
In view of Lemma \ref{Le.sec:App.44-1}
$$
\E_\zs{\p,S}\,M^2_n
\le
 \frac{\sigma_\zs{1,*}(b-a)}{d}
 \sum_{j=1}^{d}\,
\theta^{2}_\zs{j,d}=
 \frac{\sigma_\zs{1,*}(b-a)}{d}
 \Vert S\Vert^{2}_\zs{d}
\le
\frac{\sigma_\zs{1,*}(b-a)^{2}}{d}
\,,
$$
where  $\sigma_\zs{1,*}$ is given in
\eqref{sec:bound-sig-25}. Moreover, using that $\E_\zs{\p,S}\,M_n=0$, we get
$$
|\E_\zs{\p,S}\,M_\zs{n}\,\Chi_\zs{\Gamma}|\,=\,|\E_\zs{\p,S}\,M_n\,\Chi_\zs{\Gamma^c}|\,
\le 
(b-a)
\sqrt{\frac{\sigma_\zs{1,*} \P_\zs{\p,S}(\Gamma^{c})}{d}}\,.
$$
Therefore, Proposition \ref{Pr.sec:Prs.stp.times.22}
yields
\begin{equation}\label{U.12}
\lim_\zs{n\to\infty} \,n^{2\k/(2\k+1)}\,
\sup_\zs{S\in \Theta_\zs{\varepsilon,L}}\,
|\E_\zs{\p,S}\,M_n\,\Chi_\zs{\Gamma}|=0\,.
\end{equation}

\noindent 
Now, the property
\eqref{sec:Ms.6-lmt-01}
Proposition \ref{Pr.sec:App.44}
and
Lemma \ref{Le.sec:App.44-323} imply the inequality \eqref{U.2}.
Hence Theorem~\ref{Th.Up.1}.
\endproof


\noindent It is clear, that
 Theorem~\ref{Th.sec:OrIn.1} and Theorem~\ref{Th.Up.1}
  imply Theorem~\ref{Th.2.2-2}.

\bigskip

{\bf Acknowledgements.}
The last author was  partially supported 
by RFBR and CNRS (research project number 20 - 51 - 15001).

\setcounter{section}{0}
\renewcommand{\thesection}{\Alph{section}}

\section{Appendix}\label{sec:A}

\subsection{Properties of the prior distribution \ref{sec:Lo.5}}

In this section we study  properties of the distribution used in
\eqref{sec:Bayes-risk-1}.

\begin{lemma}\label{Le.sec:App.PriDens-1}
For any $\N> 2$ there exists  a continuously differentiable probability density $\rho_\zs{\N}(\cdot)$
on $\bbr$
with the support on the interval $[-\N\,,\,\N]$, i.e. $\rho_\zs{\N}(z)>0$ for $-\N<z<\N$ and   
$\rho_\zs{\N}(z)=0$ for $\vert z\vert\ge \N$,
such that for any $\N>2$ the integral
$\int^{\N}_\zs{-\N}z\rho_\zs{\N}(z)\d z=0$ and, moreover, 
$
\int^{\N}_\zs{-\N}z^{2}\rho_\zs{\N}(z)\d z
\to
1$
and
$J_\zs{\N}=\int_\zs{\bbr}\,(\dot{\rho}_\zs{\N}(z))^{2}/\rho_\zs{\N}(z)\,
\d z
\to 1$
as $\N\to\infty$.
\end{lemma}

\proof
First  we set
$
V(z)=\left(\int^{1}_\zs{-1}\,e^{-\frac{1}{1-t^{2}}}\d t\right)^{-1}\,e^{-\frac{1}{1-z^{2}}}\,\Chi_\zs{\{\vert z\vert\le 1\}}$.
\noindent
It is clear that this function is infinitely times continuously differentiable, such that
$V(z)>0$ for $\vert z\vert < 1$,
 $V(z)=0$ for $\vert z\vert \ge 1$ and $\int^{1}_\zs{-1} V(z)\d z=1$.
Now for  $\N\ge 2$ we set
$
\chi_\zs{\N}(z)=
\int^{1}_\zs{-1}\,\Chi_\zs{\{\vert z+u\vert\le N-1\}}\,V(u)\d u
=
\int_\zs{\bbr}\,\Chi_\zs{\{\vert t\vert\le N-1\}}\,V(t-z)\d t$.
Using here the properties of the function $V$ we can obtain directly that 
$\chi_\zs{\N}(z)=\chi_\zs{\N}(-z)$ for $z\in\bbr$,
$\chi_\zs{\N}(z)=1$ for $\vert z\vert\le \N-2$, $\chi_\zs{\N}(z)>0$ for $\N-2<\vert z\vert< \N$ and 
$\chi_\zs{\N}(z)=0$ for $\vert z\vert\ge \N$. Moreover, it is clear that the derivative
$
\dot{\chi}_\zs{\N}(z)=-
\int_\zs{\bbr}\,\Chi_\zs{\{\vert t\vert\le N-1\}}\,\dot{V}(t-z)\d t
=-
\int^{1}_\zs{-1}\,\Chi_\zs{\{\vert u+z\vert\le N-1\}}\,\dot{V}(u)\d u
$.
Note here that
$
\vert\dot{V}(z)\vert
\le \c_\zs{*} \sqrt{V(z)}$ for some
$\c_\zs{*}>0$.
Now through the
 Bunyakovsky - Cauchy  - Schwartz inequality we get that 
 $\dot{\chi}^{2}_\zs{\N}(z)\le 2\c_\zs{*}\chi_\zs{\N}(z)
 \quad\mbox{for}\quad
 \vert z\vert <\N$.
Now  we set
$
\rho_\zs{\N}(z)=
\left(
\int^{N}_\zs{-\N}\varphi(t)\chi_\zs{\N}(t)\d t \right)^{-1}
\varphi(z)\chi_\zs{\N}(z)$, where $\varphi(z)$  the $(0,1)$ Gaussian density,
It is clear that $\rho_\zs{\N}(z)$ is the the continuously differentiable probability density with the support $[-\N\,,\,\N]$
such that for any $\N$ the integral
$\int^{\N}_\zs{-\N}z\rho_\zs{\N}(z)\d z=0$
and  $
\int^{N}_\zs{-\N}\varphi(t)\chi_\zs{\N}(t)\d t\to 1$,
$
\int^{\N}_\zs{-\N}z^{2}\rho_\zs{\N}(z)\d z\to 1$
for $\N\to\infty$.
Moreover, the Fisher  information can be represented as 
$$
J_\zs{N}=
\left(
\int^{\N}_\zs{-\N}
\varphi(t)\chi_\zs{\N}(t)\d t \right)^{-1}
\left(
\int_\zs{\bbr}
\frac{\dot{\varphi}^{2}(z)}{\varphi(z)}\chi_\zs{\N}(z)\d z
+
\Delta_\zs{\N}
\right)
\,,
$$
where,  taking into account that  $\dot{\chi}_\zs{\N}(z)=0$ for $\vert z\vert\le \N-2$, 
$$
\Delta_\zs{\N}=
2\int_\zs{\vert z\vert \ge \N-2}\dot{\varphi}(z)\dot{\chi}_\zs{\N}(z)\d z
+
\int_\zs{\vert z\vert \ge \N-2}\,\varphi(z)\frac{\dot{\chi}^{2}_\zs{\N}(z)}{\chi_\zs{\N}(z)}\d z
\,.
$$
Therefore, $\Delta_\zs{\N}\to 0$ as $\N\to\infty$.
Hence Lemma \ref{Le.sec:App.PriDens-1}.
\endproof

\begin{lemma}\label{Le.sec:App.PriDstr-00}
The term \eqref{Psi-Def-1} is such that
$
\lim_\zs{n\to\infty}
\max_\zs{1\le j\le d}\,
\vert
\overline{\Psi}_\zs{n,j}
-
1
\vert
=0$.
\end{lemma}
\proof
First, note that 
$
y_\zs{l}=y_\zs{0}\,\prod^{l}_\zs{i=1}\,S_\zs{\kappa}(x_\zs{i})
+\sum^{l}_\zs{\iota=1}\,\prod^{l}_\zs{i=\iota+1}\,S_\zs{\kappa}(x_\zs{i})\xi_\zs{\iota}$ for $l\ge 1$.
Therefore,
$
\wt{\E}_\zs{0}\,y^{2}_\zs{l}=
y^{2}_\zs{0}\,\wt{\E}_\zs{0}\prod^{l}_\zs{i=1}\,S^{2}_\zs{\kappa}(x_\zs{i})
+\sum^{l-1}_\zs{\iota=1}\,\wt{\E}_\zs{0}\,\prod^{l}_\zs{i=\iota+1}\,S^{2}_\zs{\kappa}(x_\zs{i})+1$
and due to  \eqref{UnifUpperBound-S} we obtain that  
for any $n\ge 1$ for which $\delta^{*}_\zs{n}<1$
$$
\sup_\zs{l\ge 1}
\vert \wt{\E}_\zs{0}\,y^{2}_\zs{l}-1\vert
\le 
(\delta^{*}_\zs{n})^{2}
\frac{\left(
 y^{2}_\zs{0}+
 1
\right)
}{1-(\delta^{*}_\zs{n})^{2}}
\to 0
\quad\mbox{as}\quad
n\to\infty\,.
$$

\noindent
Since, 
$(b-a) \sum^{n}_\zs{l=1}\phi^{2}_\zs{j}(x_\zs{l})=n$,
we get  Lemma 
\ref{Le.sec:App.PriDstr-00}.
\endproof

\begin{lemma}\label{Le.sec:App.3+1}
The term $\R_\zs{0,n}$ in \eqref{sec:Lo.12} is such that
$\lim_\zs{n\to\infty}\,n^{\b}\,\R_\zs{0,n}=0$
for any $\b>0$ and $0<\rho<1$.
\end{lemma}
\proof
First note, that taking into account in the definition of term $\R_\zs{0,n}$ 
in \eqref{sec:Lo.12}, that $\vert \eta^{*}_\zs{j}\vert \le \ln n$, 
we get that
$
\R_\zs{0,n}\le 
\left(
\r+
\ln^{2} n\,
\sum^{d}_\zs{j=1}\,s^{2}_\zs{j}
\right)
\,
\mu_\zs{\kappa}(\D^{c}_\zs{n})$.
Therefore, to show this lemma it suffices to check that 
$\lim_\zs{n\to\infty}\,n^{\b}\,\mu_\zs{\kappa}(\D^{c}_\zs{n})=0$ for any $\b>0$.
To do this note, that the definition of $\D_\zs{n}$ in
\eqref{LB-1-BR-01} implies
$
\mu_\zs{\kappa}(\D^{c}_\zs{n})\le \P(\zeta_\zs{n}>\r)$
and
 $\zeta_\zs{n}=\sum^{d_\zs{n}}_\zs{j=1}\,\a_\zs{j}\kappa^{2}_\zs{j}$. So, it suffices to show that
\begin{equation}
\label{lim_prpy_1}
\lim_\zs{n\to\infty}\,
n^{\b}\,
\P(\zeta_\zs{n}>\r)
=0
\quad\mbox{for any}\quad
\b>0
\,.
\end{equation}
Indeed, first note, that  the definition \eqref{sec:Lo.5} through Lemma \ref{Le.sec:App.PriDens-1}
and the property
\eqref{sec:Lo.PR-d}
 imply
 directly
$$
\lim_\zs{n\to \infty}\,
\E\,\zeta_\zs{n}=
\lim_\zs{n\to\infty} \sum^{d}_\zs{j=1}\a_\zs{j}s^{2}_\zs{j}\E\,(\eta^{*}_\zs{1})^{2}=
\lim_\zs{n\to\infty} \sum^{d}_\zs{j=1}\a_\zs{j}s^{2}_\zs{j}=
\rho\r\,.
$$
\noindent
Setting now
$
\wt{\zeta}_\zs{n}=\zeta_\zs{n}-\E\,\zeta_\zs{n}
=(b-a)\sum^{d_\zs{n}}_\zs{j=1}\,s^{*}_\zs{j}\,\a_\zs{j}\wt{\eta}_\zs{j}/n$
and
$\wt{\eta}_\zs{j}=(\eta^{*}_\zs{j})^{2}-
\E
(\eta^{*}_\zs{j})^{2}$,
 we get  for large $n$ that 
$
\left\{
\zeta_\zs{n}
>\r
\right\}
\subset
\left\{
\wt{\zeta}_\zs{n}>
\r_\zs{1}
\right\}
$ 
for $\r_\zs{1}=\r(1-\rho)/2$.
Now the correlation inequality from
\cite{GaltchoukPergamenshchikov2013}
and the bound $\vert \wt{\eta}_\zs{j}\vert\le 2\ln^{2} n$ imply
 that for any $p\ge 2$ there exists some constant $C_\zs{p}>0$  for which
$$
\E\,\wt{\zeta}^{p}_\zs{n}\le C_\zs{p}
\frac{(\ln n)^{2p}}{n^{p}}\,
\left( \sum^{d}_\zs{j=1}\,
(s^{*}_\zs{j})^{2}\,\a^{2}_\zs{j}
\right)^{p/2}
\le C_\zs{p} \,n^{-\frac{p}{4\k+2}}\,(\ln n)^{2p}
\,,
$$
i.e.  the expectation
$\E\,\wt{\zeta}^{p}_\zs{n}\to 0$ as $n\to\infty$ and,  therefore, 
$
n^{\b}\P(\wt{\zeta}_\zs{n}>\r_\zs{1})\to 0$
as
$n\to\infty$ for any $\b>0$.
This implies \eqref{lim_prpy_1} and,  hence  Lemma \ref{Le.sec:App.3+1}. 
\endproof

\subsection{Properties of the trigonometric basis.}\label{subsec:A.2}

First we need the following lemma from 
\cite{KonevPergamenshchikov2015}.

\begin{lemma}\label{Le.sec:A.1-1}
Let $f$ be an absolutely continuous function, $f: [a,b]\to\bbr,$ with
$\|\dot{f}\|<\infty$ and $g$ be a piecewise constant function $ [a,b]\to\bbr$ of a form
$
g(x)=\sum_\zs{j=1}^{d}\,c_\zs{j}\,\chi_\zs{(z_\zs{j-1}, z_\zs{j}]}(x)$
where $c_\zs{j}$ are some constants. Then for any $\varepsilon>0,$  the function $\Delta=f-g$
satisfies the following inequalities
$$
 \|\Delta\|^{2}_\zs{d}\le (1+\wt{\varepsilon})\|\Delta\|^{2}
+
(1+\wt{\varepsilon}^{-1})\frac{(b-a)^{2}\|\dot{f}\|^{2}}{d^{2}}
\,.
$$
\end{lemma}

\begin{lemma}\label{Le.sec:A.3-00}
For any  $1\le j\le d$ 
the trigonometric Fourier  coefficients $(\theta_\zs{j,d})_\zs{1\le j\le d}$
for the functions $S$ from
the class $\cW_\zs{\k,\r}$ with $\k\ge 1$ satisfy, for any $\wt{\varepsilon}>0$,
the following inequality
$
\theta^{2}_\zs{j,d}
\,
\le\,(1+\wt{\varepsilon})
\,\theta^{2}_\zs{j}
\,
+4 \r(1+\wt{\varepsilon}^{-1})(b-a)^{2\k} d^{-2k}$.
\end{lemma}
\proof
First we represent the function $S$ as
$
S(x)=\sum^{d}_\zs{l=1}\,\theta_\zs{l}\,\phi_\zs{l}(x)
+\Delta_\zs{d}(x)$ and
$\Delta_\zs{d}(x)=\sum_\zs{l>d}\,\theta_\zs{l}\,\phi_\zs{l}(x)$, i.e.
$
\theta_\zs{j,d}=(S,\phi_\zs{j})_\zs{d}
=\theta_\zs{j}
+
(\Delta_\zs{d},\phi_\zs{j})_\zs{d}$
and, therefore,  $\forall \wt{\varepsilon}>0$ we get
$
\theta^{2}_\zs{j,d}
\le 
(1+\wt{\varepsilon})
\theta^{2}_\zs{j}
+
(1+\wt{\varepsilon}^{-1})
\Vert\Delta_\zs{d}\Vert^{2}_\zs{d}$.
 Lemma~\ref{Le.sec:A.1-1} with $g=0$ implies
$
\Vert\Delta_\zs{d}\Vert^{2}_\zs{d}\le 
2\Vert\Delta_\zs{d}\Vert^{2}
+
2(b-a)^{2} d^{-2}\Vert\dot{\Delta}_\zs{d}\Vert^{2}$.
Using here that 
$
2\pi[l/2]\ge l$ for $l\ge 2$,
we obtain that
$
\Vert\Delta_\zs{d}\Vert^{2}=\sum_\zs{l>d}\theta^{2}_\zs{l}
\le \r/\a_\zs{d} \le (b-a)^{2\k}\r d^{-2\k} $
and
\begin{equation}
\label{rest-101}
\Vert\dot{\Delta}_\zs{d}\Vert^{2}=
(2\pi/(b-a))^{2}
\sum_\zs{l>d}\,\theta^{2}_\zs{l}\,[l/2]^{2}
\le (b-a)^{2(\k-1)} d^{-2(\k-1)}\,.
\end{equation}
Hence Lemma~\ref{Le.sec:A.3-00} \endproof

\begin{lemma}\label{Le.sec:A.3}
For any $d\ge 2$ and $1\le N\le d$  
the coefficients $(\theta_\zs{j,d})_\zs{1\le j\le d}$
of functions $S$ from
the class $\cW_\zs{\r,\k}$ with $\k\ge 1$ satisfy, for any $\wt{\varepsilon}>0$,
the  inequality
$
\sum^{d}_\zs{j=N}
\theta^{2}_\zs{j,d}
\,
\le\,(1+\wt{\varepsilon})
\,\sum_\zs{j\ge N}\,\theta^{2}_\zs{j}
\,
+(1+\wt{\varepsilon}^{-1})(b-a)^{2\k}\r\,d^{-2} N^{-2(\k-1)}$.
\end{lemma}
\proof
Note that
$
 \sum^{d}_\zs{j=N}
\theta^{2}_\zs{j,d}
=\min_\zs{x_\zs{1},\ldots,x_\zs{N-1}}
\,\|S-\sum^{N-1}_\zs{j=1}\,x_\zs{j}\phi_\zs{j}
\|^{2}_\zs{d}
\le \|\Delta_\zs{N}\|^{2}_\zs{d}$
and $\Delta_\zs{N}(t)=\sum_\zs{j\ge N}\,\theta_\zs{j}\phi_\zs{j}(t)$.
 Lemma~\ref{Le.sec:A.1-1} and
\eqref{rest-101},
 imply  Lemma~\ref{Le.sec:A.3} \endproof

\subsection{Technical lemmas}
\label{subsec:App-2}

\begin{lemma}\label{Le.sec:App.44-1}
For any non random coefficients $(\u_\zs{j,l})_\zs{1\le j\le d}$
$$
\E\left(
\sum^{d}_\zs{j=1}
\u_\zs{j,l} \eta_\zs{j,d}
 \right)^{2}
 \le \sigma_\zs{1,*}
 \sum^{d}_\zs{j=1}\,
\u_\zs{j,l}^{2} 
\,,
$$
where the coefficient $\sigma_\zs{1,*}$ is given in \eqref{sec:bound-sig-25}.
\end{lemma}

\proof
Using the definition of $\eta_\zs{j,d}$ in  \eqref{sec:Ms.5} and the bounds \eqref{sec:bound-sig-25},
we get
\begin{align*}
\E
\left(
\sum^{d}_\zs{j=1}
\u_\zs{j,l} \eta_\zs{j,d}
 \right)^{2}
 &=
\frac{b-a}{d}\E
\sum^{d}_\zs{l=1}
\sigma^{2}_\zs{l}
\left(
\sum^{d}_\zs{j=1}
\u_\zs{j,l}
\phi_\zs{j}(z_\zs{l})
\right)^{2}\\[2mm]
&\le 
\sigma_\zs{1,*}
\frac{b-a}{d}
\sum^{d}_\zs{l=1}
\left(
\sum^{d}_\zs{j=1}
\u_\zs{j,l}
\phi_\zs{j}(z_\zs{l})
\right)^{2}
\,.
\end{align*}
Now, the orthonormality property \eqref{sec:Ms.4nn_1}
implies this lemma.
\endproof

\begin{lemma}\label{Le.sec:App.44-323}
For the sequence \eqref{U.7-010} the following limit property holds true
$$
\lim_\zs{n\to\infty}n^{2\k/(2\k+1)}\sup_\zs{S\in \cW_\zs{\k,\r}}\,\vert \E_\zs{\p,S}\,\U_\zs{n}\Chi_\zs{\Gamma}\vert=0
\,.
$$
\end{lemma}

\proof
First of all, note that, using the definition of $\s_\zs{j,d}$ in \eqref{sec:Ms.11}, we obtain
$$
\E_\zs{\p,S}\,\eta^2_\zs{j,d}
=
\E_\zs{\p,S}\,\s_\zs{j,d}
=
\frac{1}{d}
\sum^{d}_\zs{l=1}\,
\E_\zs{\p,S}\,
\frac{1}{H_\zs{l}}
+
\frac{1}{d}
\E_\zs{\p,S}\overline{\s}_\zs{j,d}\,,
$$
where 
$
\overline{\s}_\zs{j,d}=
\sum^{d}_\zs{l=1}\,
\sigma^2_\zs{l}\overline{\phi}_\zs{j}(x_\zs{l})$ and
$\overline{\phi}_j(z)=(b-a)\phi^2_\zs{j}(z)-1$.
Therefore, we can represent the expectation of $\U_\zs{n}$ as
$$
\E_\zs{\p,S}\,\U_\zs{n}=\frac{\Vert \wt{\lambda}\Vert^{2}}{d^{2}}\E_\zs{\p,S}\U_\zs{1,n}+
\frac{b-a}{d^{2}}
\E_\zs{\p,S}\U_\zs{2,n}\,,
$$
where $\Vert \wt{\lambda}\Vert^{2}=\sum^{d}_\zs{j=1} \wt{\lambda}^{2}(j)$,
$$
\U_\zs{1,n}=
\frac{b-a}{d}
\sum^{d}_\zs{l=1}\,
\frac{d}{H_\zs{l}}
-\varsigma_\zs{*}
\quad\mbox{and}\quad
\U_\zs{2,n}=
\sum^{d}_\zs{j=1} \wt{\lambda}^{2}(j)
\overline{\s}_\zs{j,d}\,.
$$

\noindent
Note now, that using
Proposition \ref{sec:Prs.stp.times.1}
and the dominated convergence theorem
in
the definition
\eqref{sec:Sp.8}
we obtain that
$$
\lim_\zs{n\to\infty}
\max_\zs{1\le l\le d}\,
\sup_\zs{S\in \Theta_\zs{\varepsilon,L}}\,
\sup_\zs{\p\in \cP}
\E_\zs{\p,S}\,\left\vert 
\frac{d}{H_\zs{l}}
 -(1-S^{2}(z_\zs{l}))
 \right\vert
 =0
 \,.
$$
Taking into account that for the functions from the class
\eqref{sec:Sp.1} their derivatives are uniformly bounded, we can  deduce that
$$
\lim_\zs{n\to\infty}
\sup_\zs{S\in \Theta_\zs{\varepsilon,L}}\,
\left\vert
\frac{b-a}{d}
\sum^{d}_\zs{l=1}\,
(1-S^{2}(z_\zs{l}))
-
\varsigma_\zs{*}
\right\vert
=0\,,
$$
i.e.
$
\lim_\zs{n\to\infty}
\sup_\zs{S\in \Theta_\zs{\varepsilon,L}}\,
\sup_\zs{\p\in \cP}\,
\vert
\E_\zs{\p,S}\,
\U_\zs{1,n}
\vert
=0$.
Therefore, taking into account that
$$
\limsup_\zs{n\to\infty}\,n^{2\k/(2\k+1)}\sup_\zs{S\in \Theta_\zs{\varepsilon,L}}\,
\frac{\Vert \wt{\lambda}\Vert^{2}}{d^{2}}
<\infty\,,
$$
we obtain that
\begin{equation}
\label{upper-bound-121-00}
\overline{\lim}_\zs{n\to\infty}\,n^{2\k/(2\k+1)}
\sup_\zs{S\in \Theta_\zs{\varepsilon,L}}\,
\sup_\zs{\p\in \cP}\,\vert \E_\zs{\p,S}\,\U_\zs{n}\vert
\le 
(b-a)
\overline{\lim}_\zs{n\to\infty}\,
\U^{*}_\zs{2,n}\,,
\end{equation}
where
$$
\U^{*}_\zs{2,n}
=
\frac{n^{2\k/(2\k+1)}}{d^{2}}
\sup_\zs{S\in \Theta_\zs{\varepsilon,L}}\,
\sup_\zs{\p\in \cP}\,
\vert
\E_\zs{\p,S}\U_\zs{2,n}
\vert\,.
$$

Now, using  Lemma A.2 from \cite{GaltchoukPergamenshchikov2009a} we obtain that
\begin{align*}
\left|\,
\E_\zs{\p,S}\U_\zs{2,n}
\right|\,&=
\left|\sum_{l=1}^{d}\,\E_\zs{\p,S}\,\sigma^2_l\,
\sum_{j=1}^{d}\,\wt{\lambda}^2(j)\,\overline{\phi}_\zs{j}(z_\zs{l})\right|\\[2mm]
&\le\,
d\,\sigma_\zs{1,*}\,(2^{2\k+1}+2^{\k+2}+1)
\le\,
5\,d\,\sigma_\zs{1,*}\,2^{2\k}
\,.
\end{align*}
\noindent
The definition of $\sigma_\zs{1,*}$ in
\eqref{sec:bound-sig-25} implies
$
\limsup_\zs{n\to\infty}\,d\,\sigma_\zs{1,*}<\infty$,
i.e.
$$
\limsup_\zs{n\to\infty}\,\sup_\zs{S\in \Theta_\zs{\varepsilon,L}}\,\sup_\zs{\p\in \cP}
\left|\,
\E_\zs{\p,S}\U_\zs{2,n}
\right|\,
<\infty\,.
$$
\noindent
Therefore, the using this bound in \eqref{upper-bound-121-00}
implies
$$
\overline{\lim}_\zs{n\to\infty}\,n^{2\k/(2\k+1)}
\sup_\zs{S\in \Theta_\zs{\varepsilon,L}}\,
\sup_\zs{\p\in \cP}\,\vert \E_\zs{\p,S}\,\U_\zs{n}\vert
=0\,.
$$
Using the inequality (A.4) from \cite{ArkounBruaPergamenchtchikov2019}
we get
$
\E_\zs{\p,S}
\eta^4_\zs{j,d}\le 64 \v^{*} \sigma^{2}_\zs{1,*}$,
where the coefficient $\v^{*}$ is given in \eqref{sec:def-eta-prs-25-03}. From this we obtain, that
\begin{align*}
\E_\zs{\p,S}|\U_\zs{n}|\Chi_\zs{\Gamma^c}&\le 
\frac{(b-a)}{d}\sum^{d}_\zs{j=1}\E_\zs{\p,S}\eta^2_\zs{j,d}\Chi_\zs{\Gamma^c}+
\varsigma_\zs{*}\P_\zs{\p,S}(\Gamma^c)\\
&\le \frac{8 \sigma_\zs{1,*}(b-a)\sqrt{\v^{*}}}{d}
\sqrt{\P_\zs{\p,S}(\Gamma^c)}+\varsigma_\zs{*}
\P_\zs{\p,S}(\Gamma^c)\,.
\end{align*}
\noindent
So,  Proposition \ref{Pr.sec:Prs.stp.times.22}
implies
$
\lim_\zs{n\to\infty}\,n^{2\k/(2\k+1)}
\sup_\zs{S\in \cW_\zs{\k,\r}} \E_\zs{\p,S}|\U_\zs{n}|\Chi_\zs{\Gamma^c}
=0$.
Hence
Lemma \ref{Le.sec:App.44-323}.
\endproof

\medskip

\bibliographystyle{plain}


\end{document}